\documentclass[showpacs,preprintnumbers,amsmath,amssymb,12pt]{article}
\usepackage{graphicx,amsfonts,amssymb,amsthm,amsmath,graphics,epsfig,pstricks,fancyhdr,fancybox}
\usepackage{dcolumn}
\usepackage{bm}

\newtheorem{thm}{Teorema}[section]

\newtheorem{prop}[thm]{Proposition}

\newcommand{\refeq}[1]{~(\ref{#1})}
\newcommand{\myref}[1]{~\ref{#1}}
\newcommand{\mycite}[1]{~\cite{#1}}

\newcommand{\RE}{\bm{R}}
\newcommand{\rv}{\textit{rv}}
\newcommand{\id}{\textit{id}}
\newcommand{\iid}{\textit{iid}}
\newcommand{\sd}{\textit{sd}}

\newcommand{\pdf}{\textit{pdf}}
\newcommand{\cdf}{\textit{cdf}}

\newcommand{\chf}{\textit{chf}}

\newcommand{\Pqo}{\bm{P}\hbox{-\emph{a.s.}}}
\newcommand{\eqd}{\stackrel{d}{=}}

\newcommand{\PR}[1]{\bm{P}\left\{{#1}\right\}}
\newcommand{\EXP}[1]{\bm{E}\left[{#1}\right]}
\newcommand{\VAR}[1]{\bm{V}\left[{#1}\right]}
\newcommand{\cov}[2]{\bm{cov}\left[{#1},{#2}\right]}

\newcommand{\poiss}{\mathfrak{P}}

\newcommand{\erl}{\mathfrak{E}}

\newcommand{\bin}{\mathfrak{B}}

\newcommand{\indim}{\noindent{\bf Proof:}\hspace{0.2cm}}
\newcommand{\findim}{\hfill$\square$\vspace{0.3cm}\noindent}

\newcommand{\heavi}{\vartheta}

\def\ito{It\={o}}

\begin{document}
\thispagestyle{empty}

\title{\Huge \textbf{Correlated Poisson processes \\ and self-decomposable laws}
\author{Nicola \textsc{Cufaro Petroni}\footnote{cufaro@ba.infn.it}  \\
Dipartimento di \textsl{Matematica} and \textsl{TIRES}, Universit\`a di Bari\\
\textsl{INFN} Sezione di Bari\\ \vspace{7pt}
via E. Orabona 4, 70125 Bari, Italy\\
Piergiacomo \textsc{Sabino}\footnote{piergiacomo.sabino@uniper.energy}\\
%\textsl{RQPR} Quantitative Risk Modelling and Analytics \\
\textsl{Uniper} Global Commodities SE\\
\vspace{5pt}
 Holzstrasse 6, 40221 D\"usseldorf, Germany
} }

\date{}

\maketitle

\begin{abstract}
\noindent We analyze a method to produce pairs of non independent
Poisson processes $M(t),N(t)$ from positively correlated,
self-de\-com\-po\-sa\-ble, exponential renewals. In particular the
present paper provides the family of copulas pairing the renewals,
along with the closed form for the joint distribution $p_{m,n}(s,t)$
of the pair $\big(M(s),N(t)\big)$, an outcome which turns out to be
instrumental to produce explicit algorithms for applications in
finance and queuing theory. We finally discuss the cross-correlation
properties of the two processes and the relative timing of their
jumps
\end{abstract}

%\vspace{5pt}

%\noindent \emph{PACS}:

%\vspace{5pt}

%\noindent \emph{Key words}: Geometric Wiener process

%\newpage

%\begin{center}

%\end{center}

%\newpage

%\tableofcontents

%\newpage

\section{Introduction}

Recent studies have shown that the spot dynamics of commodity
markets displays mean reversion, seasonality and
jumps\mycite{CarteaFigueroa}, and some methodologies have been
proposed to take dependency into account based on correlation and
co-integration\mycite{DH13}. However, these approaches can become
mathematically cumbersome and non-treatable when leaving the
Gaussian-{\ito} world. In this context it has been indeed recently
proposed\mycite{energyconf} to consider $2$-dimensional jump
diffusion processes with a $2$-dimensional Gaussian and a
$2$-dimensional compound Poisson component, and, as also suggested
in different circumstances\mycite{iyer}, we show here that a
revealing approach to model the dependency of the $2$-dimensional
Poisson processes can be supplied on the ground of the
\emph{self-decomposability} of the exponential random variables used
for its construction.

The present paper is in particular devoted to find both the copula
function pairing our correlated renewals, and an explicit form for
the joint distribution $p_{m,n}(s,t)$ of our pair of correlated
Poisson processes $M(s),N(t)$: this will prove to be instrumental to
produce the efficient algorithms that can be used in financial
applications\mycite{energyconf}. If indeed the pairs of correlated,
exponential random variables (\rv's) $(X_k,Y_k)$ -- used to produce
the renewals in our processes -- are interpreted as random waiting
times with \emph{random delays}, the proposed model can help
describing their co-movement and can answer some common questions
arising in the financial context:
    \begin{itemize}
        \item Once a financial institution defaults how long should one wait for a dependent institution to default too?
        \item A market receives a news interpreted as a shock: how long should one wait to see the propagation of that shock onto a dependent market?
        \item What is the impact of the correlations among the shocks for different insurance companies on a fair assessment of the risk of losses\mycite{meyers}?
    \end{itemize}
It is worth noticing, moreover, that we achieve our aim of producing
a $2$-dimensional Poisson process with \emph{dependent} marginals
without resorting to an \emph{a priori} copula (distributional)
approach: the dependence among arrival times will indeed be made
explicit in terms of combinations of \rv's, and we only recover and
discuss the corresponding copula functions as an outcome of this
model. As a consequence, because of this $\Pqo$ relationship between
the random times, the two Poisson processes can be seen as linked
with a form of \emph{co-integration} between their jumps. Similar
models -- albeit rather less sophisticated -- were also used in
order to model a multi-component reliability system\mycite{iyer},
while the so-called \emph{Common Poisson Shock
Models}\mycite{LindskogMcNeil} are in fact quite different from that
presented here

%In particular we can find close-form formulas for vanilla options
%hence the price and the Greeks of spread options can be calculated
%in close form using the Margrabe formula [5] (if the strike is zero)
%or some well known approximations as in Deng et al. [8]. In any case
%our approach implies an explicit algorithm for the simulation of the
%dependent Poisson processes and can be used in Monte Carlo
%simulations.

The main practical consequence of our results is then that the price
and the Greeks of the spread options considered in the
applications\mycite{energyconf} can be calculated in \emph{closed
form} using either the Margrabe formula (if the strike is zero), or
some well known approximation\mycite{deng}. In any case our model
entails explicit algorithms for the simulation of correlated Poisson
processes, and can be used in the \emph{Monte Carlo} simulations. An
extension to the multi-dimensional case, as well as to different
dynamics other than Poisson, will be considered in future studies;
but, under the assumption that only two underlyings have jump
component, the price and the Greeks of spread options can be
obtained even now by the moment-matching methodology recently
proposed in\mycite{PellegrinoSabino}

The paper is organized as follows: in the Section\myref{sdlaws} we
first show how (hitherto \emph{positively}) correlated exponential
\rv's can be deduced from the self-decomposability of their laws;
then in the Section\myref{copfun} we briefly discuss the copula
functions produced by this model. By using pairs of these
exponential \rv's as correlated renewals, in the Section\myref{pcpp}
we subsequently produce a 2-dimensional Poisson process with
correlated components, and in the Section\myref{joint} we explicitly
deduce their joint distributions. Finally in the
Section\myref{discuss} the cross-correlation properties of the
Poisson processes are briefly analyzed, and the relative timing of
their jumps is used to shed new light on the dependence mechanism of
a model allowing for the possibility of a delayed propagation of
correlated shocks. We conclude by pointing out first that we would
also be able to produce other correlated \rv's (Erlang, Gamma,
\emph{EPT}...) by making use, once more, of their
self-decomposability; and then that the results of this paper should
also be extended to \emph{negatively} correlated renewals, a
possibility -- not open to other procedures -- that will be
postponed to future inquiries. Lengthy proofs are confined in the
Appendices, together with a few details about the notation adopted
throughout the paper

\section{Correlation from self-decomposability}\label{sdlaws}

\subsection{Joint distributions}

A law with density (\pdf) $f(x)$ and characteristic function (\chf)
$\varphi(u)$ is said to be \emph{self-decomposable}
(\sd)\mycite{sato,cufaro08} when for every $0<a<1$ we can find
another law with \pdf\ $g_a(x)$ and \chf\ $\chi_a(u)$ such that
\begin{equation*}
    \varphi(u)=\varphi(au)\chi_a(u)
\end{equation*}
This is
%not at all a trivial requirement because, albeit for every
%$\varphi(u)$ and $0<a<1$ it would be always possible to take
%\begin{equation*}
%    \chi_a(u)=\frac{\varphi(u)}{\varphi(au)}
%\end{equation*}
%it is by no means apparent that the right-hand side of this relation
%would always be a good \chf: this happens only under special
%circumstances which give meaning to the definition of
%self-decomposability and select
a well known family of laws with many relevant properties. We will
also say that a random variable (\rv) $X$ with \pdf $f(x)$ and \chf\
$\varphi(u)$ is \sd\ when its law is \sd: looking at the definition
this means that for every $0<a<1$ we can always find two
\emph{independent} \rv's $Y$ (with the same law of $X$) and $Z_a$
(with \pdf\ $g_a(x)$ and \chf\ $\chi_a(u)$) such that
\begin{equation*}
    X\eqd aY+Z_a
\end{equation*}
We can look at this, however, also from a different perspective: if
$Y$ is \sd, and to the extent that, for $0<a<1$, an independent
$Z_a$ with the suitable law is known, we can define a third \rv
\begin{equation*}
   X\equiv aY+Z_a\qquad\quad\Pqo
\end{equation*}
being sure that it will have the same law as $Y$. In the following
we will mainly adopt this second standpoint

We turn now, for later convenience, to give the joint laws of the
triplet $(X,Y,Z_a)$: for the \chf\ $\psi(u,v,w)$ we easily find from
the independence of $Y$ and $Z_a$ that
\begin{eqnarray*}
  \psi(u,v,w)& = &\EXP{e^{i(uX+vY+wZ_a)}}\\
  &=&\varphi(au+v)\chi_a(u+w)=\varphi(au+v)\frac{\varphi(u+w)}{\varphi(a(u+w))}
\end{eqnarray*}
while the marginal, joint \chf's of the pairs $(X,Y)$ and  $(X,Z_a)$
respectively are
\begin{eqnarray*}
    \phi(u,v)&=&\psi(u,v,0)=\varphi(au+v)\frac{\varphi(u)}{\varphi(au)}\\
    \omega(u,w)&=&\psi(u,0,w)=\varphi(au)\frac{\varphi(u+w)}{\varphi(a(u+w))}
\end{eqnarray*}
As for the \pdf\ $\kappa(x,y,z)$ on the other hand , by taking
$s=au+v$, we have (we neglect the integration limits whenever they
extend to all $\RE$)
\begin{eqnarray*}
  \kappa(x,y,z) &=& \frac{1}{(2\pi)^3}\int \!\!du\int \!\!dv \int \!\!dw \,e^{-i(ux+vy+wz)}\varphi(au+v)\chi_a(u+w)\\
   &=&\frac{1}{(2\pi)^2}\int \!\!du\int \!\!dw
   \,e^{-i(ux+wz)}\chi_a(u+w)  \frac{1}{2\pi}\int
   e^{-ivy}\varphi(au+v)\,dv \\
   &=&\frac{1}{(2\pi)^2}\int \!\!du\int \!\!dw
   \,e^{-i(ux+wz)}\chi_a(u+w)e^{-iyau}  \frac{1}{2\pi}\int
   e^{-isy}\varphi(s)\,ds \\
   &=& f(y)\frac{1}{(2\pi)^2}\int \!\!du\int \!\!dw
   \,e^{-i[u(x-ay)+wz]}\chi_a(u+w)
\end{eqnarray*}
and again with $s=u+w,\;t=w$
\begin{eqnarray*}
  \kappa(x,y,z) &=& f(y)\frac{1}{(2\pi)^2}\int \!\!ds\int \!\!dt
   \,e^{-i[(x-ay)(s-t)+zt]}\chi_a(s) \\
   &=&f(y)\,\frac{1}{2\pi}\int \!\!ds
   \,e^{-is(x-ay)}\chi_a(s)\frac{1}{2\pi}\int
   \!\!dt\,e^{-it[z-(x-ay)]} \\
   &=&f(y)\,g_a(x-ay)\,\delta[z-(x-ay)]
\end{eqnarray*}
so that the marginal, joint \pdf's of $(X,Y)$ and $(X,Z_a)$ will
respectively be
\begin{equation}\label{XYjoint}
    h(x,y)=f(y)\,g_a(x-ay)\qquad\qquad
    \ell(x,z)=\frac{1}{a}\,f\!\left(\frac{x-z}{a}\right)g_a(z)
\end{equation}
Finally the joint cumulative distribution function (\cdf) of $(X,Y)$
is
\begin{equation*}
  H(x,y)=\int_{-\infty}^yf(y')\,G_a(x-ay')\,dy'\qquad\qquad  G_a(z)=\int_{-\infty}^zg_a(z')\,dz'
\end{equation*}
where $G_a(z)$ is the  \cdf\ of $Z_a$. The particular form of
$H(x,y)$ will be instrumental in finding the copula
functions\mycite{nels} eventually pairing $X$ and $Y$

We can finally also calculate the correlation coefficients $r_{XY}$
and
 $r_{XZ_a}$: if we put $\EXP{X}=\EXP{Y}=\mu$ and $\VAR{X}=\VAR{Y}=\sigma^2$,
from the $Y,Z_a$ independence we have
\begin{equation*}
    \EXP{XY}=\EXP{(aY+Z_a)Y}=a\sigma^2+\mu^2
\end{equation*}
and hence $r_{XY}=a$. In a similar vein, to calculate $r_{XZ_a}$ we
first remark that $\VAR{X}=a^2\VAR{Y}+\VAR{Z_a}$, namely
$\VAR{Z_a}=(1-a^2)\sigma^2$, and then from
\begin{equation*}
  \EXP{XZ_a} = \EXP{(aY+Z_a)Z_a}=(1-a^2)\sigma^2+(1-a)\mu^2
\end{equation*}
we finally find $r_{XZ_a}=1-a^2$

\subsection{An example: the exponential laws
$\erl_1(\lambda)$}\label{example}

It is well known that the exponential laws $\erl_1(\lambda)$ with
\pdf\ and \chf\ (see Appendix\myref{notations} for the notations
adopted from now on)
\begin{equation*}
    \lambda f_1(\lambda x)=\lambda e^{-\lambda x}\heavi(x) \qquad\qquad \varphi_1(u/\lambda)=\frac{\lambda}{\lambda-iu}
\end{equation*}
are a typical example of \sd\ laws\mycite{sato}, and in this case we
can explicitly give the law of $Z_a$: we have indeed
\begin{equation}\label{chichf}
    \chi_a(u)=\frac{\varphi_1(u/\lambda)}{\varphi_1(au/\lambda)}=\frac{\lambda-iau}{\lambda-iu}=a+(1-a)\frac{\lambda}{\lambda-iu}=a+(1-a)\varphi_1(u/\lambda)
\end{equation}
which (for $0<a<1$) is a mixture of a law $\bm\delta_0$ degenerate
in $0$, and an exponential $\erl_1(\lambda)$, namely
\begin{equation*}
    Z_a\sim a\bm\delta_0+(1-a)\erl_1(\lambda)
\end{equation*}
so that its \pdf\ and \cdf\ respectively are
\begin{eqnarray*}
    g_a(z)&=&a\delta(z)+(1-a)\lambda e^{-\lambda z}\heavi(z)\\
    G_a(z)&=&\left[a+(1-a)(1-e^{-\lambda z})\right]\heavi(z)
\end{eqnarray*}
It is also easy to prove on the other hand that this coincides with
the law of the product of two other independent \rv's: an
exponential $Z\sim\erl_1(\lambda)$, and a Bernoulli
$B(1)\sim\bin(1,1-a)$ with $a=\PR{B(1)=0}$, so that we can always
write
\begin{equation*}
    Z_a=B(1)\,Z
\end{equation*}
In short, given two exponential \rv's $Y\sim\erl_1(\lambda)$ and
$Z\sim\erl_1(\lambda)$, and a Bernoulli $B(1)\sim\bin(1,1-a)$, all
three mutually independent, the \rv\ $X$ defined as
\begin{equation}\label{expdecomp}
    X\equiv aY+B(1)Z
\end{equation}
is again an exponential $\erl_1(\lambda)$. From\refeq{XYjoint} we
also find that the joint \pdf\ of $X,Y$ is
\begin{equation*}
    h(x,y)=\lambda e^{-\lambda y}\heavi(y)\left[a\delta(x-ay)+(1-a)\lambda e^{-\lambda(x-ay)}\heavi(x-ay)\right]
\end{equation*}
and hence its joint \cdf\ is
\begin{eqnarray*}
  H(x,y) &=& \int_{-\infty}^y\lambda e^{-\lambda y'}\heavi(y')\left[a+(1-a)(1-e^{-\lambda(x-ay')})\right]\heavi(x-ay')\,dy' \\
   &=&\heavi\left(y\wedge\frac{x}{a}\right)\bigg[\left(1- e^{-\lambda\left(y\wedge\frac{x}{a}\right)}\right)-e^{-\lambda x}\left(1- e^{-\lambda(1-a)\left(y\wedge\frac{x}{a}\right)}\right)\bigg]
\end{eqnarray*}
Of course this is far from the only possible joint law with
exponential marginals (see also Section\myref{cop}), but it is
noticeable because it traces its origins back to a model of
self-decomposability of the exponentials. As for the correlations
among $X,Y$ and $Z$, we already know that $r_{XY}=a$. For $r_{XZ}$
we first find that
\begin{equation*}
    \EXP{XZ}=\EXP{(aY+Z_a)Z}=\frac{2-a}{\lambda^2}
\end{equation*}
and then that
\begin{equation*}
    r_{XZ}=1-a=1-r_{XY}
\end{equation*}
so that for our three exponentials in\refeq{expdecomp} we eventually
have
\begin{equation*}
   r_{XY}+r_{XZ}=1\qquad\quad r_{XY}=a\qquad\quad r_{YZ}=0
\end{equation*}

\subsection{Positively correlated exponential \rv's}

It is apparent now from the discussion in the previous section that
the self-de\-com\-po\-sa\-bi\-li\-ty of the exponential laws
$\erl_1(\lambda)$ can be turned into a simple procedure to generate
identically distributed and correlated \rv's: given
$Y\sim\erl_1(\lambda)$, in order to produce another
$X\sim\erl_1(\lambda)$ with correlation $0<a<1$, it would be enough
to take $Z\sim\erl_1(\lambda)$ and $B(1)\sim\bin(1,1-a)$ independent
from $Y$ and define $X$ as in\refeq{expdecomp}. In other words $X$
will be nothing else than the exponential $Y$ \emph{down}
$a$-\emph{rescaled}, plus another independent exponential $Z$
\emph{randomly intermittent with frequency} $1-a$. The
self-decomposability of the exponential laws ensures then that, for
every $0<a<1$, also $X$ marginally is an $\erl_1(\lambda)$. Remark
that we would not have the same result by taking more naive
combinations of $Y$ and $Z$. Consider for instance the sum
$aY+(1-a)Z$ of our two independent, exponential \rv's: in this case,
since $aY\sim\erl_1\left(\frac{\lambda}{a}\right)$ and
$(1-a)Z\sim\erl_1\left(\frac{\lambda}{1-a}\right)$, the law of
$aY+(1-a)Z$ would be
$\erl_1\left(\frac{\lambda}{a}\right)*\erl_1\left(\frac{\lambda}{1-a}\right)$,
which is neither an exponential $\erl_1(\lambda)$, nor even in
general an Erlang $\erl_2$ because of the difference between the two
parameters
%: it in fact goes into an $\erl_2(2\lambda)$ only in the
%limit $a\to\! ^1/_2$
\begin{figure}
\begin{center}
\includegraphics*[width=14cm]{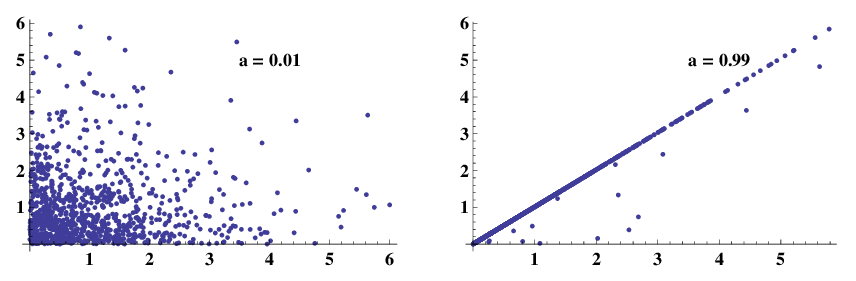}
\caption{The pairs $(X_k,Y_k)$ with correlation $0.01$ and
$0.99$}\label{renew}
\end{center}
\end{figure}

The proposed procedure can now be adapted to generate a sequence of
independent pairs of exponential \rv's, with correlated components,
$(X_k,Y_k)$, $k=1,2,\ldots$ (or, if we prefer, $X_k,Z_k$) that will
act in the subsequent sections as renewals for a two-dimensional
point (Poisson) process: take indeed $0<a<1$, produce two
independent \id\ exponentials $Y,Z$ and another independent
Bernoulli $B(1)$, then define $X=aY+B(1)Z$ and take the pair $X,Y$.
By independently replicating this procedure we will get a sequence
of \iid\ two-dimensional pairs $(X_k,Y_k)$ that will be used later
to generate a two-dimensional Poisson process with correlated
renewals. Of course -- not surprisingly -- the case of uncorrelated
pairs of renewals $(X_k,Y_k)$, and hence of independent Poisson
processes, is retrieved from our model in the limit $a\to0$, because
in this case we just have $X=Z$ which is by definition independent
from $Y$. On the other hand it is also apparent that in the opposite
limit $a\to1$ (namely when $r_{XY}$ goes to 1) we tend to have
$X=Y,\;\Pqo$ so that the time pairs will fall precisely on the
diagonal of the two-dimensional time, and the two Poisson processes
will simply $\Pqo$ coincide. To see it from another standpoint we
could look to some simulation of the pairs $(X_k,Y_k)$: for small
correlations $a$ the scatter-plot of our pairs $(X_k,Y_k)$ tends to
evenly spread out within the first quadrant without any apparent
hint to some forme of dependence; on the other hand for $a$ near to
1 the points tend to cluster together along the diagonal, as can be
seen in the Figure\myref{renew}

\section{Copulas for bivariate exponentials}\label{copfun}

\subsection{A family of copula functions}\label{cop}

From the discussion in the Section\myref{example} we know that the
pair $X,Y$ of correlated \rv's deduced from their
self-decomposability has the joint \cdf\
\begin{equation}\label{jointcdf}
    H(x,y)=\heavi\left(y\wedge\frac{x}{a}\right)\bigg[\left(1- e^{-\lambda\left(y\wedge\frac{x}{a}\right)}\right)
        -e^{-\lambda x}\left(1- e^{-\lambda(1-a)\left(y\wedge\frac{x}{a}\right)}\right)\bigg]
\end{equation}
with the exponential marginal \cdf's (the notation is here slightly
simplified)
\begin{equation}\label{expcdf}
    F(x)=\heavi(x)\left(1-e^{-\lambda x}\right)\qquad\qquad G(y)=\heavi(y)\left(1-e^{-\lambda y}\right)
\end{equation}
To find out the \emph{copula function} $C(u,v)$\mycite{nels} pairing
$X,Y$ we then first remark that
\begin{equation*}
    e^{-\lambda x}=1-F(x)\qquad\quad e^{-\lambda
    y}=1-G(y)\qquad\quad e^{-a\lambda
y}=\left[1-G(y)\right]^a
\end{equation*}
while
\begin{eqnarray*}
  0\le\frac{x}{a}\le y\quad &\Longleftrightarrow&\quad e^{-a\lambda y}\le e^{-\lambda x}\quad\; \Longleftrightarrow\quad\; \left[1-G(y)\right]^a\le1-F(x) \\
  0\le y\le\frac{x}{a}\quad &\Longleftrightarrow&\quad e^{-a\lambda y}\ge e^{-\lambda x}\quad\; \Longleftrightarrow\quad\; \left[1-G(y)\right]^a\ge1-F(x)
\end{eqnarray*}
and then that our joint \cdf\refeq{jointcdf} takes the
form\begin{eqnarray*}
  H(x,y)=F(x)\qquad\qquad\qquad\qquad\qquad\qquad\qquad &\hbox{for}&\!\!\!\! \left[1-G(y)\right]^a\le1-F(x) \\
  H(x,y)=F(x)-[1-G(y)]\left(1-\frac{1-F(x)}{[1-G(y)]^a}\right) &\hbox{for}&\!\!\!\! \left[1-G(y)\right]^a\ge1-F(x)
\end{eqnarray*}
which can also be conveniently summarized as
\begin{equation*}
    H(x,y)=F(x)-[1-G(y)]\left(1-\frac{1-F(x)}{[1-G(y)]^a}\right)^+
\end{equation*}
As a consequence we get the following family of copula functions
\begin{equation}\label{copula}
    C_a(u,v)=u-(1-v)\left[1-\frac{1-u}{(1-v)^a}\right]^+=u-\frac{[(1-v)^a-(1-u)]^+}{(1-v)^{a-1}}
\end{equation}
which for $0\le a\le1$ runs between two extremal copulas
\begin{eqnarray*}
  C_0(u,v) &=& uv\qquad\qquad\quad\, \hbox{independent marginals}\\
  C_1(u,v) &=& u\wedge v\qquad\qquad \hbox{fully positively correlated marginals}
\end{eqnarray*}
It is easy to see that $C_1(u,v)$ also coincides with the
Fr\échet-H\"offding upper bound $\overline{C}(u,v)$ for copulas (see
Section\myref{frechet})

\subsection{Bivariate exponential
distributions}\label{alternative}

Several examples -- all different from\refeq{copula} -- of bivariate
distributions with exponential marginals $\erl_1(\lambda)$ and
$\erl_1(\mu)$ can be found in the literature\mycite{nels,balak}.
First we find the \emph{Gumbel bivariate exponential
distribution}\mycite{nels,gumb} with $0\le a\le1$ and
\begin{eqnarray*}
    H(x,y)&=&\heavi(x)\heavi(y)\left(1-e^{-\lambda x}-e^{-\mu y}+e^{-\lambda x-\mu
    y-a\lambda\mu\,xy}\right)\\
    h(x,y)&=&\heavi(x)\heavi(y)\left[a(\lambda x+\mu
    y+a\lambda\mu\,xy)+1-a\right]e^{-\lambda x-\mu
    y-a\lambda\mu\,xy} \\
    C_a(u,v)&=&u+v-1+(1-u)(1-v)e^{-\frac{a}{\lambda\mu}\ln(1-u)\ln(1-v)}
\end{eqnarray*}
It is apparent that $C_0(u,v)=uv$ gives the independent
exponentials, while
\begin{equation*}
    C_1(u,v)=u+v-1+(1-u)(1-v)e^{-\frac{1}{\lambda\mu}\ln(1-u)\ln(1-v)}
\end{equation*}
does not seem to correspond to some notable copula. Then there is
the \emph{Marshall-Olkin bivariate exponential
distribution}\mycite{nels,marsh} with $0\le a,b\le1$ and
\begin{eqnarray*}
    H(x,y)\!\!\!\!&=&\!\!\!\!\left\{
               \begin{array}{ll}
                 \!\!\heavi(x)\heavi(y)(1-e^{-\lambda x})^{1-a}(1-e^{-\mu y}) & \hbox{if}\;(1-e^{-\lambda x})^a\ge(1-e^{-\mu y})^b\\
                 \!\!\heavi(x)\heavi(y)(1-e^{-\lambda x})(1-e^{-\mu y})^{1-b} & \hbox{if}\;(1-e^{-\lambda x})^a\le(1-e^{-\mu y})^b
               \end{array}
             \right.
    \\
    h(x,y)\!\!\!\!&=&\!\!\!\!\left\{
               \begin{array}{ll}
                 \!\!\frac{1-a}{(1-e^{-\lambda x})^a}\,\heavi(x)\lambda e^{-\lambda x}\,\heavi(y)\mu e^{-\mu y} & \hbox{if}\;(1-e^{-\lambda x})^a\ge(1-e^{-\mu y})^b\\
                 \!\!\frac{1-b}{(1-e^{-\mu y})^b}\,\heavi(x)\lambda e^{-\lambda x}\,\heavi(y)\mu e^{-\mu y} & \hbox{if}\;(1-e^{-\lambda x})^a\le(1-e^{-\mu y})^b
               \end{array}
             \right. \\
   \lefteqn{ C_{a,b}(u,v)=(u^{1-a}v)\wedge(uv^{1-b})\left\{
                                                \begin{array}{ll}
                                                  u^{1-a}v & \quad\hbox{when}\;\; u^a\ge v^b \\
                                                  uv^{1-b} & \quad\hbox{when}\;\; u^a\le v^b
                                                \end{array}
                                              \right.}
\end{eqnarray*}
In this case $C_{0,0}(u,v)=uv$ again is the independent copula,
while $C_{1,1}(u,v)=u\wedge v$ is the Fr\échet-H\"offding upper
bound $\overline{C}(u,v)$ (see Section\myref{frechet}): apart from
these extremal values, however, also this Marshall-Olkin copula
differs from\refeq{copula}. A third family of copulas can finally be
traced back to the \emph{Raftery bivariate exponential
distribution}\mycite{nels,raft}: in this case the copula functions
are
\begin{equation*}
    C_a(u,v)=u\wedge v+\frac{a}{2-a}(uv)^{\frac{1}{a}}\left[1-(u\vee v)^{1-\frac{2}{a}}\right]
\end{equation*}
and correspond to the case of correlated exponential \rv's $X,Y$
which are produced by three independent exponential \rv's $U,V$ and
$Z$ according to the definitions
\begin{equation*}
    X\equiv aU+B(1)Z\qquad\quad Y\equiv aV+B(1)Z
\end{equation*}
Here, at variance with our model based on self-decomposability, the
correlation is apparently produced by the presence of the same
exponential \rv\ $Z$ in both the right-hand sides of the
definitions. In short, it results from these examples that our
family of copulas\refeq{copula} seems not to have been used in
advance to couple pairs of marginal exponentials

\subsection{Fr\échet-H\"offding bounds}\label{frechet}

It is well known known\mycite{nels} that every copula function
$C(u,v)$ falls between the Fr\échet-H\"offding bounds
\begin{equation*}
   \underline{C}(u,v)= (u+v-1)^+\le C(u,v)\le u\wedge v=\overline{C}(u,v)
\end{equation*}
and we have also found in the Section\myref{cop} that the copula
$C_1(u,v)$ for our fully correlated ($r_{XY}=1$) exponential
marginals coincides with the Fr\échet-H\"offding upper bound. By
keeping in mind a possible generalization of our model to the case
of \emph{negatively correlated} exponentials, we will briefly recall
in this section a few general features of the joint \cdf's
$H(x,y)=C(F(x),G(y))$ produced by the pairing of two given \cdf's
$F(x)$ and $G(x)$ by means of the Fr\échet-H\"offding lower and
upper bounds

Let us suppose for simplicity that $F(x)$ and $G(x)$ are continuous
and strictly increasing functions so that the inverse functions
exist, and consider first the lower bound copula
$\underline{C}(u,v)= (u+v-1)^+$: in that case the condition
$F(x)+G(y)\ge1$ is equivalent to both the inequalities
\begin{equation*}
    x\ge\beta(y)=F^{-1}(1-G(y))\qquad\qquad y\ge\alpha(x)=G^{-1}(1-F(x))
\end{equation*}
and hence from $H(x,y)=(F(x)+G(y)-1)^+$ we first have
\begin{eqnarray*}
  \partial_xH(x,y) &=& \left\{
        \begin{array}{ll}
          f(x) &\qquad \hbox{if $\;\;y\ge\alpha(x)$} \\
          0 & \qquad\hbox{if $\;\;y<\alpha(x)$}
        \end{array}
      \right. \\
  \partial_yH(x,y) &=& \left\{
        \begin{array}{ll}
          g(y) &\qquad \hbox{if $\;\;x\ge\beta(y)$} \\
          0 & \qquad\hbox{if $\;\;x<\beta(y)$}
        \end{array}
      \right.
\end{eqnarray*}
where $f(x)$ and $g(y)$ are the corresponding marginal \pdf's, and
then the joint \pdf\ is
\begin{equation}\label{jointcorr}
    h(x,y)=\partial_x\partial_yH(x,y)=f(x)\delta(y-\alpha(x))=g(y)\delta(x-\beta(y))
\end{equation}
As a consequence we can say that the joint laws produced by the
copula $\underline{C}(u,v)$ describe pairs of coupled \rv's $X,Y$
satisfying $\Pqo$ the functional relations
\begin{equation}\label{functrel}
    X=\beta(Y)=F^{-1}(1-G(Y))\qquad\qquad Y=\alpha(X)=G^{-1}(1-F(X))
\end{equation}
A formally identical result can be proved for the case of the upper
bound 2-copula $\overline{C}(u,v)=u\wedge v$ but for the fact that
now the functions $\alpha(x)$ and $\beta(y)$ must be redefined as
\begin{equation*}
    \alpha(x)=G^{-1}(F(x))\qquad\qquad \beta(y)=F^{-1}(G(y))
\end{equation*}

In the case of the lower bound copula $\underline{C}(u,v)$ it is
interesting to remark now that when for instance $F(x)$ and $G(y)$
are Gaussian \cdf's the functions $\alpha(x)$ and $\beta(y)$ are
linear with negative proportionality coefficients, so that the pair
$X,Y$ is perfectly anti-correlated with $r_{XY}=-1$. The same
happens also in the case of a pair of Student laws of the same
order. This is true indeed for every other pair of marginal laws of
the \emph{same type} and with \emph{support coincident with} $\RE$.
On the other hand when the marginals either are not of the same
type, or have an unbounded support strictly contained in $\RE$ (as
happens for exponential laws), they apparently can not reciprocally
be in a linear relation with negative proportionality coefficient,
and hence can not be totally \emph{linearly} anti-correlated. In
this case it can still be proved by means of H\"offding's Lemma
(see\mycite{nels} p.\ 190) that the minimal correlation is reached
by means of the lower bound copula $\underline{C}$, but now
$\alpha(x)$ and $\beta(y)$ can no longer be \emph{linear} functions,
and $r_{XY}$ will be strictly larger than $-1$. By taking indeed the
Fr\échet-H\"offding lower bound $\underline{C}(u,v)=(u+v-1)^+$ as
the copula for our exponentials\refeq{expcdf} we would find the
\pdf\refeq{jointcorr} and the functional relations\refeq{functrel}
where now
\begin{equation*}
    \alpha(x)=-\frac{1}{\lambda}\ln\left(1-e^{-\lambda x}\right)\qquad\quad\beta(y)=-\frac{1}{\lambda}\ln\left(1-e^{-\lambda y}\right)
\end{equation*}
and a short calculation would then show that in this case
\begin{equation*}
    r_{XY}=1-\frac{\pi^2}{6}\approx-0.645
\end{equation*}
so that this minimal anti-correlation allowed for exponential \rv's
would in any case be larger than $-1$. It could in fact be proved in
general (see\mycite{nels} p.\ 30-32) that, when $X$ and $Y$ are
continuous, $Y$ is almost surely a decreasing function of $X$ if and
only if the copula of $X$ and $Y$ is $\underline{C}$. Random
variables with copula $\underline{C}$ are often called
countermonotonic. We postpone to a subsequent enquiry a detailed
study of negatively correlated exponentials

\section{Correlated Poisson processes}\label{pcpp}

Following the discussion of Section\myref{sdlaws} it is easy now to
produce a sequence of \rv's by independently iterating the
definition\refeq{expdecomp}
\begin{equation}\label{pairs}
    X_k=aY_k+B_k(1)Z_k\qquad\quad k=1,2,\ldots
\end{equation}
in such a way that for every $k$: $X_k,Y_k,Z_k$ are
$\erl_1(\lambda)$, $B_k(1)$ are $\bin(1,1-a)$, and $Y_k,Z_k,B_k(1)$
are mutually independent. The pairs $(X_k,Y_k)$ instead will be
$a$-cor\-re\-la\-ted for every $k$. Add moreover
$X_0=Y_0=Z_0=0,\;\Pqo$ to the list, and take then the point
processes for $n=0,1,2,\ldots$
\begin{equation}\label{pointproc}
  T_n = \sum_{k=0}^nX_k\qquad\quad S_n = \frac{\lambda}{\mu}\sum_{k=0}^nY_k\qquad\quad
R_n = \sum_{k=0}^nZ_k
\end{equation}
Since the $X_k\sim\erl_1(\lambda)$ are \iid\ \rv's we know that
$T_n\sim\erl_n(\lambda)$ are distributed as Erlang (gamma) laws with
\pdf's $\lambda f_n(\lambda x)$ and \chf's $\varphi_n(u/\lambda)$
(see Appendix\myref{notations} for notations) where it is understood
that $T_0\sim\erl_0=\bm\delta_0$. In a similar way we can argue that
$S_n\sim\erl_n(\mu)$ and $R_n\sim\erl_n(\lambda)$. We will finally
denote with $N(t)\sim\poiss(\lambda t)$ and $M(t)\sim\poiss(\mu t)$
the \emph{correlated Poisson processes} associated respectively to
$T_n$ and $S_n$
\begin{figure}
\begin{center}
\includegraphics*[width=14cm]{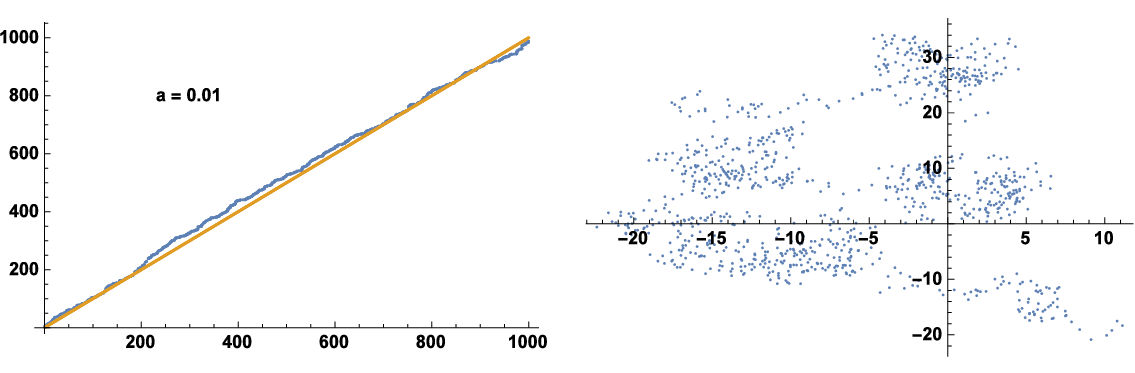}
\caption{Sample pairs of the two-dimensional point process
$(T_n,S_n)$ with correlation $r_{XY}=0.01$: on the left the points
are compared with the average trend
$\big(\frac{n}{\lambda},\frac{n}{\mu}\big)$; on the right they are
instead plotted after centering around these averages}\label{corr01}
\end{center}
\end{figure}

In order to get a first look to these processes we generate
$n=1\,000$ pairs $(X_k,Y_k)$ with the associated two dimensional
point process $(T_n,S_n)$, and then we simulate the corresponding
Poisson processes $N(t)$ and $M(t)$. The pairs $(T_n,S_n)$ are first
plotted along with their average time increases
$\big(\frac{n}{\lambda},\frac{n}{\mu}\big)$, and then after
centering around these averages, namely as
\begin{equation*}
    T_n-\frac{n}{\lambda}\qquad\quad S_n-\frac{n}{\mu}\qquad\quad
    n=1,2,\ldots, 1\,000
\end{equation*}
In this second rendering the random behavior is magnified by
consistently reducing the plot scale to a suitable size. In the same
way for the Poisson processes we first show samples of the pair
$N(t),M(t)$, and then that of their \emph{compensated} versions
$\widetilde{N}(t)=N(t)-\lambda t$ and $\widetilde{M}(t)=M(t)-\mu t$

\begin{figure}
\begin{center}
\includegraphics*[width=14cm]{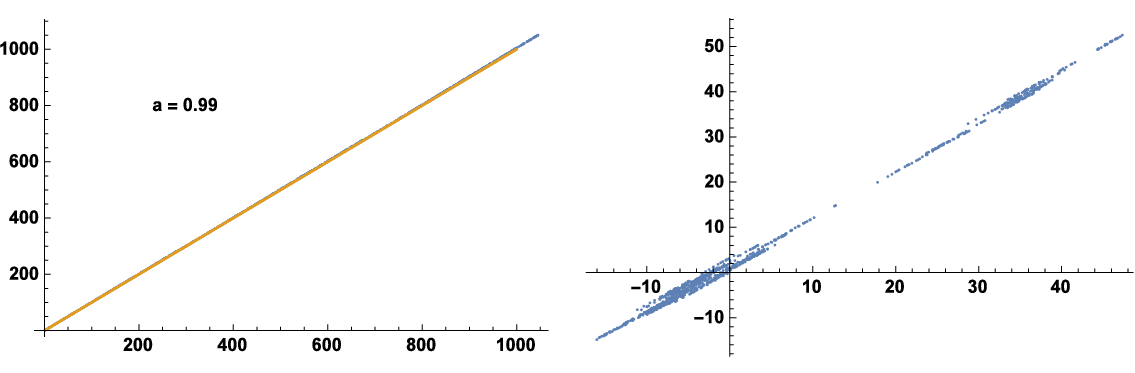}
\caption{Sample pairs of the two-dimensional point process
$(T_n,S_n)$ with correlation $r_{XY}=0.99$: on the left the points
are compared with the average trend
$\big(\frac{n}{\lambda},\frac{n}{\mu}\big)$; on the right they are
instead plotted after centering around these averages}\label{corr99}
\end{center}
\end{figure}

\begin{figure}
\begin{center}
\includegraphics*[width=14cm]{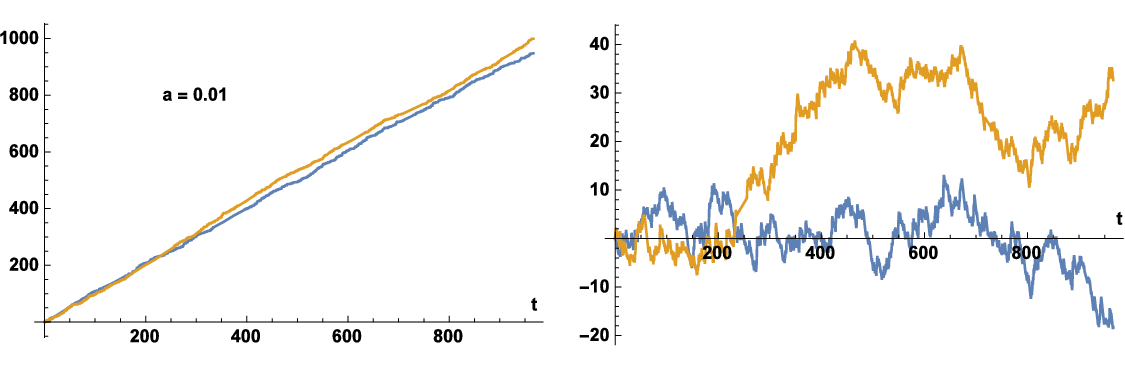}
\caption{On the left sample paths are shown of the two Poisson
processes $N(t)$ and $M(t)$ with correlation $r_{XY}=0.01$; on the
right we instead have the corresponding compensated Poisson
processes $\widetilde{N}(t)$ and $\widetilde{M}(t)$}\label{Pcorr01}
\end{center}
\end{figure}

\begin{figure}
\begin{center}
\includegraphics*[width=14cm]{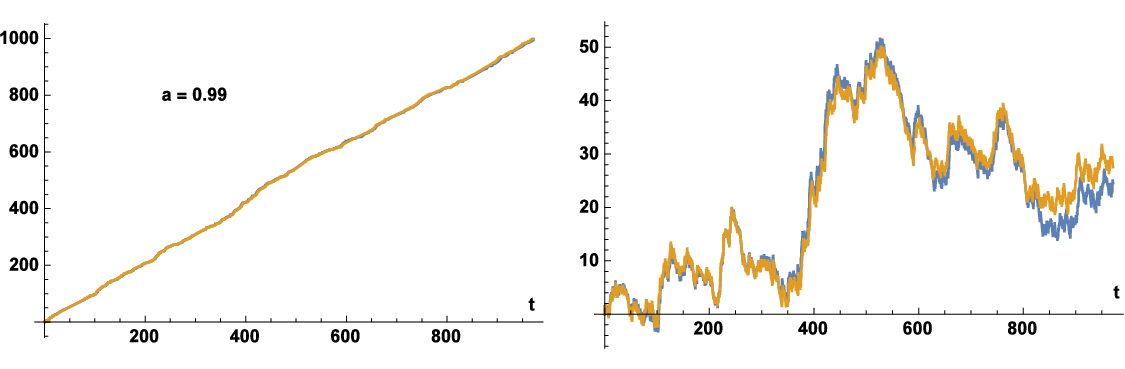}
\caption{On the left sample paths are shown of the two Poisson
processes $N(t)$ and $M(t)$ with correlation $r_{XY}=0.99$; on the
right we instead have the corresponding compensated Poisson
processes $\widetilde{N}(t)$ and $\widetilde{M}(t)$}\label{Pcorr99}
\end{center}
\end{figure}

In the Figure\myref{corr01} we plotted the two-dimensional point
process $(T_n,S_n)$ with $\lambda=\mu=1$ and $a=r_{XY}=0.01$: since
the correlation among the renewals is negligible the right-hand
plots (centered around the averages) apparently show a random
behavior. In the Figure\myref{corr99} instead we took
$a=r_{XY}=0.99$, namely we generated strongly and positively
correlated renewals. In this second case, as it was to be expected,
the centered time pairs fall into line among themselves. As for the
Poisson processes themselves, in the Figure\myref{Pcorr01} the
trajectories on the left hand side have $a=0.01$ correlation and
look fairly independent, after suitable compensation, on the right
hand side. In the Figure\myref{Pcorr99} instead we took a
correlation $a=0.99$ and the compensated trajectories are now almost
superimposed. Remark as on the left-hand sides of these figures both
the Poisson processes and the time pairs appear to be quite near to
one another, and to their averages because of a scale effect which
is eliminated by compensation and centering in the corresponding
right-hand sides
\begin{prop}
The \rv's
\begin{equation*}
    \zeta_n\equiv\sum_{k=0}^nB_k(1)Z_k
\end{equation*}
turn out to be the sum of a (random) binomial number
$B(n)\sim\bin(n,1-a)$ of \iid\ exponentials $\erl_1(\lambda)$, and
hence they follow an Erlang law with a random index $B(n)$ (here
$B(0)=0$), namely
\begin{equation*}
    \zeta_n=\sum_{k=0}^{B(n)}Z_k\sim\erl_{B(n)}(\lambda)
\end{equation*}
\end{prop}
 \indim
This is better seen from the point of view of the mixtures by
remarking that, if $\varphi_1(u/\lambda)$ is the \chf\ of
$\erl_1(\lambda)$, we have from\refeq{chichf} (see also
Appendix\myref{notations})
\begin{eqnarray*}
  \varphi_{\zeta_n}(u) \!\!&=& \!\!\EXP{e^{iu\zeta_n}}=\EXP{\prod_{k=0}^ne^{iuB_k(1)Z_k}}= \prod_{k=0}^n\EXP{e^{iuB_k(1)Z_k}}\\
   &=&\!\!\left[a+(1-a)\varphi_1\left(\frac{u}{\lambda}\right)\right]^n=\sum_{k=0}^n\binom{n}{k}a^{n-k}(1-a)^k\varphi_1^k\left(\frac{u}{\lambda}\right)\\
   &=&\sum_{k=0}^n\beta_k(n)\varphi_k\left(\frac{u}{\lambda}\right)
\end{eqnarray*}
which, if $\varphi_k(u)=[\varphi_1(u)]^k$ are the \chf\ of
$\erl_k(1)$, eventually is a mixture of Erlang laws
$\erl_k(\lambda)$ with the binomial weights $\beta_k(n)$. It is
understood here that $\varphi_1^0(u)=1$, so that
$\erl_0(\lambda)=\bm\delta_0$ and $f_0(x)=\delta(x)$ (see
Appendix\myref{notations})
 \findim

\noindent A straightforward consequence of the previous proposition
(which apparently just amounts to acknowledge a subordination) is
that now from
\begin{equation*}
    \sum_{k=0}^nX_k=a\sum_{k=0}^nY_k+\sum_{k=0}^nB_k(1)Z_k
\end{equation*}
we will also have
\begin{equation}\label{TSrelation}
    T_n=\frac{a\mu}{\lambda}S_n+\zeta_n=\frac{a\mu}{\lambda}S_n+\sum_{k=0}^{B(n)}Z_k=\frac{a\mu}{\lambda}S_n+R_{B(n)}
\end{equation}
where $R_{B(n)}\sim\erl_{B(n)}(\lambda)$ is the point process $R_n$
with a random index $B(n)$.
It is worthwhile to notice that the previous results also
substantiate the well known fact that the Erlang \rv's are
self-decomposable too: the explicit knowledge of the $\zeta_n$ law
allows indeed to construct pairs of dependent Erlang \rv's with
correlation $a$

\section{The joint distribution}\label{joint}

Our main task is now to explicitly calculate the joint distribution
of our Poisson processes at arbitrary times $s,t\ge0$ and
$n,m=0,1,2,\ldots$
\begin{eqnarray*}
    p_{m,n}(s,t)&=&\PR{M(s)=m,\,N(t)=n}\\
    &=&\PR{S_m\le s<S_{m+1},\,\,T_n\le t<T_{n+1}}
\end{eqnarray*}
and to this effect we first remark (in a slightly simplified
notation) that
\begin{eqnarray}
  p_{m,n}\!\! &=&\!\! \PR{M(s)\ge m,\,N(t)\ge n} - \PR{M(s)\ge m+1,\,N(t)\ge n}\nonumber\\
   \lefteqn{-\PR{M(s)\ge m,\,N(t)\ge n+1} + \PR{M(s)\ge m+1,\,N(t)\ge
   n+1}}\nonumber\\
   &=&\!\!q_{m,n}-q_{m+1,n}-q_{m,n+1}+q_{m+1,n+1}\label{pmn}
\end{eqnarray}
where
\begin{equation*}
   q_{m,n}(s,t)=\PR{M(s)\ge m,\,N(t)\ge n}=\PR{S_m\le s,\,T_n\le t}
\end{equation*}
so that by taking
\begin{equation*}
    w=\frac{\lambda r}{a}\qquad\quad
    y=\frac{\lambda t}{a}\qquad\quad z=\frac{\lambda t-a\mu s}{a}<y
\end{equation*}
from\refeq{TSrelation} we are reduced to calculate (see also
Appendix\myref{notations})
\begin{eqnarray}
  q_{m,n} &=&\PR{S_m\le s,\,\,\frac{a\mu}{\lambda}S_n+R_{B(n)}\le t}\nonumber \\
   &=&\sum_{\ell=0}^n\beta_\ell(n)\int_0^\infty dr \lambda f_\ell(\lambda r)\nonumber\\
   &&\qquad\qquad\qquad\PR{\left.S_m\le s,\,\frac{a\mu}{\lambda} S_n+R_{B(n)}\le
   t\,\right|R_\ell=r,B(n)=\ell}\nonumber \\
   &=&\lambda\sum_{\ell=0}^n\beta_\ell(n)\int_0^t dr f_\ell(\lambda r)\PR{S_m\le
   s,\,S_n\le\lambda\frac{t-r}{a\mu}}\nonumber\\
%   &=&a\sum_{\ell=0}^n\beta_\ell(n)\int_0^y dw f_\ell(aw)\PR{S_m\le \frac{y-z}{\mu},\,S_n\le
%   \frac{y-w}{\mu}}\nonumber\\
   &=&a\int_0^y dw\, h_n(aw)\PR{S_m\le \frac{y-z}{\mu },\,S_n\le
   \frac{y-w}{\mu }}\label{qmn}
\end{eqnarray}
where $\lambda h_n(\lambda x)$ is the (Erlang binomial mixture)
\pdf\ of $R_{B(n)}$
\begin{prop}\label{prop}
For $n,m=0,1,2,\ldots$ and $\rho,\tau\ge0$ we have
\begin{align*}
%  \PR{S_m\le \rho,\,S_n\le \tau}\!\! &=&\!\! \PR{M(\rho\wedge \tau)\ge m\vee n}\\
%   &&\qquad\qquad\qquad+\left[\Theta_{n-m}\heavi(\tau-\rho)+\Theta_{m-n}\heavi(\rho-\tau)\right]\cdot \\
%   &&\sum_{k=m\wedge
%   n}^{(m\vee n)-1}\!\!\PR{M(|\rho-\tau|)\ge(m\vee n)-k}\PR{M(\rho\wedge \tau)=k}
%   \\
   &\qquad\quad\PR{S_m\le \rho,\,S_n\le \tau}\,\,=\,\,\Pi_{m\vee n}\big(\mu(\rho\wedge\tau)\big)\\
   &\!\!+\big[\Theta_{n-m}\heavi(\tau-\rho)+\Theta_{m-n}\heavi(\rho-\tau)\big]\!\!\sum_{k=m\wedge
   n}^{(m\vee n)-1}\Pi_{(m\vee
   n)-k}\big(\mu|\rho-\tau|\big)\pi_k\big(\mu(\rho\wedge\tau)\big)
\end{align*}
with the notations adopted in the Appendix\myref{notations} for the
Poisson laws
\end{prop}
 \indim
See Appendix\myref{proofprop} for a detailed proof
 \findim

\noindent Of course in\refeq{qmn} we take in particular
\begin{equation*}
    \rho=\frac{y-z}{\mu }=s\qquad\quad\tau=\frac{y-w}{\mu }=\lambda\frac{t-r}{a\mu}
\end{equation*}
It is apparent that this result will be instrumental to calculate
first $q_{m,n}(s,t)$ in\refeq{qmn}, and then the distributions
$p_{m,n}(s,t)$
\begin{prop}\label{lemma}
If $\;a\mu s \ge\lambda t$, then $p_{m,n}(s,t)=0$ whenever $m<n$
\end{prop}
 \indim
Since our renewals $X_k,Y_k,Z_k$ are all non-negative \rv's, the
point processes are always non-decreasing
\begin{equation*}
    S_m \le S_{m+1\qquad\quad }T_n \le T_{n+1}\qquad\quad
    m,n=0,1,2,\ldots
\end{equation*}
while
\begin{equation*}
    T_n=\frac{a\mu}{\lambda} S_n+\sum_{k=0}^{B(n)}Z_k\ge\frac{a\mu}{\lambda} S_n
\end{equation*}
Now, if $M(s)=m$ and $N(t)=n$, we must have both $S_m\le s<S_{m+1}$
and $T_n\le t<T_{n+1}$. Suppose now $0\le m<n$, namely $ m+1\le n$
and $S_{m+1}\le S_n$: then
\begin{equation*}
    \frac{a\mu}{\lambda} s<\frac{a\mu}{\lambda} \,S_{m+1}\le\frac{a\mu}{\lambda} \,S_n\le
    T_n\le t
\end{equation*}
which apparently contradicts the hypothesized inequality
 \findim

\noindent As a consequence when $a\mu s\ge\lambda t$ we can always
restrict our calculations to the case $m\ge n\ge 0$. We can now
finally state our complete results about the joint distributions
$p_{m,n}(s,t)$
\begin{prop}\label{case1}
Take for short
\begin{equation*}
    y=\frac{\lambda t}{a}>0\qquad\quad z=\frac{\lambda t-a\mu
    s}{a}<y
\end{equation*}
Then, when $\;a\mu s>\lambda t$, namely $z<0$, we have
\begin{equation}\label{pmn1}
  p_{m,n}(y,z)=\left\{\begin{array}{ll}
                     0 & \qquad n>m\ge0\\
                     Q_{n,n}(y,z) & \qquad m=n\ge0 \\
                     Q_{m,n}(y,z)-Q_{m,n+1}(y,z) & \qquad m>n\ge0
                   \end{array}
                 \right.
\end{equation}
where we defined
\begin{equation}\label{Q}
  Q_{m,n}(y,z)= a\int_0^ydw\,
  h_n(aw)\sum_{k=n}^{m}\pi_{m-k}(w-z)\pi_k(y-w)\qquad m\ge n\ge0
\end{equation}
When instead $\;a\mu s<\lambda t$, and hence $0<z<y$, we have
\begin{align}
  &  \qquad\quad p_{m,n}(y,z)\label{pmn2}\\
  &  =\left\{\begin{array}{ll}
                    \! A_{m,n}(y,z)-A_{m,n+1}(y,z)+B_{m,n}(y,z)-B_{m,n-1}(y,z) &
                     \quad n>m\ge0 \\
                     \! A_{n,n}(y,z)-A_{n,n+1}(y,z)+B_{n,n}(y,z)+C_{n,n}(y,z) & \quad m=n\ge0 \\
                     \! A_{m,n}(y,z)-A_{m,n+1}(y,z)+C_{m,n}(y,z)-C_{m,n+1}(y,z) & \quad m>n\ge0
                   \end{array}
                 \right.\nonumber
\end{align}
where we defined
\begin{align}
 & A_{m,n}(y,z)=
  a\!\int_0^zdw\,h_n(aw)\pi_m(y-z)\qquad\qquad\qquad\qquad\;\;
  n,m\ge0\label{A1}\\
 & B_{m,n}(y,z)=
  a\!\int_0^zdw\,h_{n+1}(aw)\sum_{k=0}^{n-m}\pi_k(z-w)\pi_m(y-z)\quad\, n\ge
  m\ge0\label{B1}\\
 & C_{m,n}(y,z)=a\int_z^ydw\,h_n(aw)\sum_{k=n}^m\pi_{m-k}(w-z)\pi_k(y-w)\quad m\ge n\ge 1\label{C1}
\end{align}
while for $m\ge n=0$ we always have $C_{m,0}(y,z)=0$. Moreover both
the results for $z<0$, and for $0<z<y$ connect with continuity in
$z=0$ in the sense that
\begin{equation*}
    p_{m,n}(y,0^-)=p_{m,n}(y,0^+)\qquad\quad m,n\ge0
\end{equation*}

\end{prop}
 \indim Take first the case $a\mu s>\lambda t$, namely $z<0$, and recall
that for the integration variable in\refeq{qmn} it is $0\le w\le y$.
As a consequence, when the Proposition\myref{prop} in used
in\refeq{qmn}, we will always have
\begin{equation}\label{ineq}
    0\le \tau=\frac{y-w}{\mu}\le\frac{y-z}{\mu}=s=\rho
\end{equation}
On the other hand, since the conditions of the
Proposition\myref{lemma} are met, we can also restrict ourselves to
evaluate $p_{m,n}(s,t)$ for $0\le n\le m$. Then, by considering
separately the cases ${m=n\ge0}$ and ${m>n\ge0}$, from\refeq{qmn}
and from the Proposition\myref{prop} we first calculate $q_{m,n},
q_{m+1,n}, q_{m,n+1}$ and $q_{m+1,n+1}$, and finally (lengthy
algebraic details can be found in the Appendix\myref{proof1})
from\refeq{pmn} we find\refeq{pmn1}.

When on the other hand $a\mu s<\lambda t$ (namely $y>z>0$ and
$0<w<y$) and we use Proposition\myref{prop} in\refeq{qmn}, instead
of\refeq{ineq} we find
\begin{equation}\label{ineq2}
   0\le \tau=\frac{y-w}{\mu}\qquad\quad0\le \rho=s=\frac{y-z}{\mu}
\end{equation}
so that $\rho$ and $\tau$ can now be in an order whatsoever. As a
consequence Proposition\myref{lemma} does not hold, and we must
consider all the possible orderings of $m,n$. Following then the
same line of reasoning as before, and always taking separately the
different $n,m$ orderings, a tedious calculation (see
Appendix\myref{proof1}) gives first the $q$'s from\refeq{qmn}, and
eventually the $p$'s of our proposition from\refeq{pmn}

We finally show that the values of $p_{m,n}(y,z)$ separately listed
in the Proposition\myref{case1} for $z<0$ and $z>0$ connect with
continuity in $z=0$, in the sense that for every $y>0$
\begin{equation*}
    p_{m,n}(y,0^-)=p_{m,n}(y,0^+)
\end{equation*}
For $z<0$ (namely $a\mu s>\lambda t$) the results are given
in\refeq{pmn1} and\refeq{Q}, so that for $z\uparrow0^-$ and every
$m,n\ge0$, we simply have
\begin{align}
 & p_{m,n}(y,0^-) = 0\qquad\qquad\qquad\qquad\qquad\qquad\qquad\qquad\quad\,\, n>m\ge0\label{pz<0n>m}\\
 & p_{n,n}(y,0^-) = a\int_0^y\!\!dw\,h_n(aw)\,\pi_0(w)\,\pi_n(y-w)\qquad\qquad\; n=m\ge0\label{pz<0n=m}\\
 & p_{m,n}(y,0^-) = a\int_0^y\!\!dw\Bigg\{h_{n+1}(aw)\pi_{m-n}(w)\pi_n(y-w)\ \quad m>n\ge0\label{pz<0n<m} \\
    & \qquad\qquad\qquad+\big[h_n(aw)-h_{n+1}(aw)\big]\sum_{k=n}^{m}\pi_{m-k}(w)\pi_k(y-w)\Bigg\}\nonumber
\end{align}
On the other hand, when $z>0$ (namely $a\mu s<\lambda t$) we
have\refeq{pmn2},\refeq{A1},\refeq{B1} and\refeq{C1}, so that now
$z$ appears also as an integration limit, and some care should be
exercised for $z\downarrow0^+$. When indeed the integrand contains
the distribution $\delta(x)$, as in fact happens in every first term
of $h_n(x)$ which is $\beta_0(n)\delta(x)=a^n\delta(x)$ (see also
Appendix\myref{notations}), we have for every regular function
$\xi(x)$
\begin{equation*}
    \lim_{z\downarrow0^+}\int_0^z\xi(x)\delta(x)\,dx=\xi(0)\qquad\qquad\lim_{z\downarrow0^+}\int_z^y\xi(x)\delta(x)\,dx=0
\end{equation*}
As a consequence we have
\begin{eqnarray*}
    \lim_{z\downarrow0^+}\int_0^zdx\,\xi(x)h_n(x)&=&a^n\xi(0)\\
    \lim_{z\downarrow0^+}\int_z^ydx\,\xi(x)h_n(x)&=&\int_0^ydx\,\xi(x)\sum_{k=1}^n\beta_k(n)f_k(x)
\end{eqnarray*}
With this provisions in mind, it is then only a question of sheer
calculation (see Appendix\myref{proof1}) to show that the
$p_{m,n}(y,z)$ for $z>0$ as given
in\refeq{pmn2},\refeq{A1},\refeq{B1} and\refeq{C1} correctly
converge to the values\refeq{pz<0n>m},\refeq{pz<0n=m}
and\refeq{pz<0n<m} for every possible ordering of $n,m$. For
instance for $n>m\ge0$, with $x=aw$ and recalling also that
$\pi_k(0)=\delta_{k,0}$ (so that $\pi_{n-m}(0)=0$ because $n>m$), in
the limit $z\downarrow0^+$ we immediately have
\begin{eqnarray*}
    \lefteqn{p_{m,n}(y,0^+)}\\
    &=&\!\!\pi_m(y)\int_0^{0^+}\!\!dx\,\left[a^n-a^{n+1}+a^n\pi_{n-m}(0)-(a^n-a^{n+1})\sum_{k=0}^{n-m}\pi_k(0)\right]\delta(x) \\
    &=&\!\!\pi_m(y)\left[a^n-a^{n+1}-(a^n-a^{n+1})\sum_{k=0}^{n-m}\delta_{k,0}\right]=0
\end{eqnarray*}
and so on for the other two cases
 \findim

%\vfill\eject

\begin{prop}\label{case2}
The terms $Q,A,B$ and $C$ in the Proposition\myref{case1} can be
expressed in terms of finite combinations of elementary functions:
when $z<0$ (namely $\;a\mu s>\lambda t$) we have for $m\ge n\ge0$
\begin{eqnarray}\label{Qmn}
  Q_{m,n}(y,z) &=&\sum_{k=n}^m\sum_{j=k}^m\frac{(-1)^{j-k}}{a^j}\binom{j}{k}\sum_{\ell=0}^n\beta_\ell(n)\\
   &&\qquad\qquad\qquad\pi_{m-j}(y-z)\pi_{j+\ell}(ay)\,\Phi(j+1;j+\ell+1;ay)\nonumber
\end{eqnarray}
Here and in the following $\Phi(\alpha;\beta;x)$ are confluent
hypergeometric functions. When instead $z>0$ (namely $\;a\mu
s<\lambda t$) we have for every $n,m\ge0$
\begin{equation}
  A_{m,n}(y,z)=\pi_m(y-z)\sum_{k=0}^n\beta_k(n)\left[1+\pi_k(az)-\sum_{j=0}^k\pi_j(az)\right]\label{A}
\end{equation}
while for $n\ge m\ge 0$ it is
\begin{eqnarray}
  B_{m,n}(y,z)&=&\pi_m(y-z)\sum_{k=0}^{n-m}\pi_k\left(z\right)\sum_{\ell=0}^{n+1}\beta_\ell(n+1)\label{B}\\
  &&\qquad\qquad\qquad\qquad\qquad \frac{(az)^\ell
  k!}{(k+\ell)!}\,\Phi\left(\ell,k+\ell+1,(1-a)z\right)\nonumber
\end{eqnarray}
and for $m\ge n\ge 1$ (for $m\ge n=0$ we always have
$C_{m,0}(y,z)=0$) it is
\begin{align}
 & C_{m,n}(y,z)=\frac{e^{-(1-a)(y-z)}}{a^m}\sum_{\ell=1}^n\beta_\ell(n)\sum_{k=n}^m\sum_{j=0}^{\ell-1}(-1)^{\ell-1-j}\binom{k+\ell-j-1}{k}\nonumber\\
  &\qquad
  \pi_j(ay)\pi_{m+\ell-j}(a(y-z))\Phi(k+\ell-j,m+\ell-j+1,a(y-z))\label{C}
\end{align}
Finally, since the parameters $\alpha,\beta$ of the
$\Phi(\alpha,\beta,x)$ involved in the previous equations are
integer numbers with $0\le \alpha<\beta$, our confluent
hypergeometric functions are just finite combinations of powers and
exponentials according to the following formulas
\begin{align*}
  & \Phi(0\,,\beta,x)=1 \qquad\qquad\qquad\qquad\qquad\qquad\qquad\qquad\qquad\qquad\quad \beta>\alpha=0\\
  & \Phi(\alpha,\beta,x)= e^x\sum_{\gamma=1}^{\alpha}(-1)^{\alpha-\gamma}\binom{\beta-\gamma-1}{\beta-\alpha-1}\frac{\pi_{\gamma-1}(x)}{\pi_{\beta-1}(x)}\,\Pi_{\alpha-\gamma+1}(x)
   \;\quad
   \beta>\alpha\ge 1
\end{align*}
\end{prop}
 \indim The detailed proof unfolds along a sequence of integrations based on known results and is here omitted for the sake of
brevity (see Appendix\myref{proof2})
 \findim

%\vfill\eject

\noindent We end this section with a short list of a few explicit
examples of joint probabilities holding in the region $a\mu
s\ge\lambda t$:
\begin{align*}
 & p_{0,0}(s,t) = e^{-\mu s} \\
 & p_{1,0}(s,t) = \frac{e^{-\mu s}}{a}\big[(1-a)(1-e^{-\lambda t})+a\mu s-\lambda t\big] \\
 & p_{1,1}(s,t) = \frac{e^{-\mu s}}{a}\big[\lambda t-(1-a)(1-e^{-\lambda t})\big] \\
 & p_{2,0}(s,t) =\frac{ e^{-\mu s}}{2a^2}\big[2(1-a)(1+a\mu s)(1-e^{-\lambda t})+(a\mu s-\lambda t)^2-2(1-a)\lambda t\big]\\
 & p_{2,1}(s,t) = \frac{e^{-\mu s}}{a^2}\big[(1-a)(a-4-(1-a)\lambda t-a\mu s)(1-e^{-\lambda t})\\
               &\qquad\qquad\qquad\qquad\qquad\qquad\qquad\qquad+\,\lambda t(a^2-5a+4+a\mu s-\lambda t)\big]\\
 & p_{2,2}(s,t) =\frac{ e^{-\mu s}}{2a^2}\big[2(1-a)(3-a+(1-a)\lambda t)(1-e^{-\lambda t})\\
               &\qquad\qquad\qquad\qquad\qquad\qquad\qquad\qquad+\,\lambda t(\lambda t-2(1-a)(3-a))\big]
\end{align*}

\section{Cross-correlations and relative timing}\label{discuss}

In this section we will briefly discuss the main correlation
properties of the two processes. We first of all look at the point
processes and we remark that, by recasting the self-decomposability
equation\refeq{pairs} in the form
\begin{equation*}
    X_k=\frac{a\mu}{\lambda}\,W_k+B_k(1)Z_k\qquad\quad
    W_k=\frac{\lambda}{\mu}\,Y_k\sim\erl_1(\mu)
\end{equation*}
the point processes appear as
\begin{equation*}
    T_n=\sum_{k=0}^nX_k\qquad\quad S_n=\sum_{k=0}^nW_k
\end{equation*}
where now $X_k\sim\erl_1(\lambda)$ and $W_k\sim\erl_1(\mu)$ play at
once the role of the correlated renewals. It is interesting to point
out then that, at variance with other models\mycite{LindskogMcNeil},
we are no longer tied to take truly \emph{coincident} shocks: we
will show indeed that with non-zero probabilities the values of the
paired, and correlated renewals $X_k,W_k$ (waiting times) can be in
an order whatsoever, and they would almost never coincide. As a
consequence the propagation of the shocks from a process to the
other will quite plausibly happen with delays whose random sizes
(and directions) could also be modeled by suitably choosing our
parameters $a,\lambda$ and $\mu$. And moreover the random times
$T_n$ and $S_n$ will be correlated by the summing up of the
renewals, but will never fall at the same instant. This
\emph{relative timing} apparently allows for an enhanced flexibility
of the model in the practical applications because we no longer have
to rely on common shocks, but rather on correlated and randomly
delayed ones
\begin{figure}
\begin{center}
\includegraphics*[width=10cm]{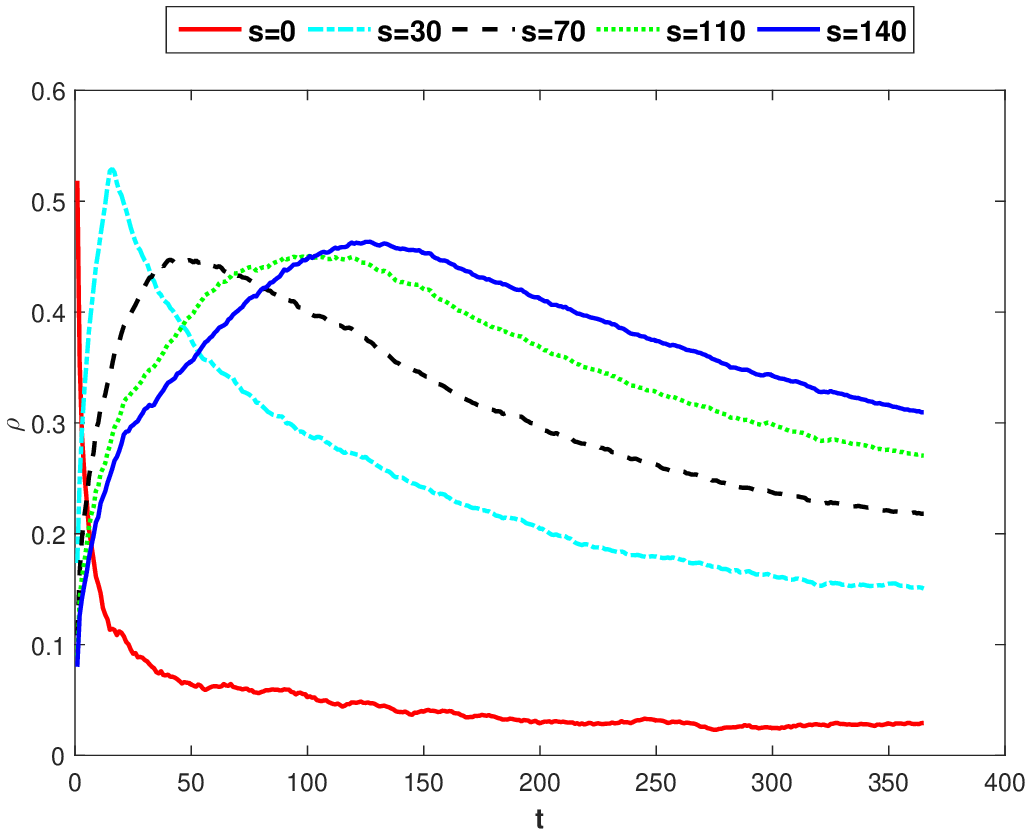}
\caption{Cross-correlation $\rho(s,t)$ of the two Poisson processes
$M(s)$ and $N(t)$ estimated by Monte Carlo simulations. Here
$a=0.5$, and $\lambda=\mu=20$, namely
$^{a\mu}/_\lambda=0.5<1$}\label{crosscorr05}
\end{center}
\end{figure}

More precisely we can now single out two possible regimes for our
processes: $^{a\mu}/_\lambda\le1$ and $^{a\mu}/_\lambda>1$. It is
then easy to see that for every $k=1,2,\ldots$
\begin{equation*}
    X_k=\frac{a\mu}{\lambda}\,W_k+B_k(1)Z_k\ge\frac{a\mu}{\lambda}\,W_k>W_k\qquad\quad\frac{a\mu}{\lambda}>1
\end{equation*}
and hence we first of all have
\begin{equation*}
    \PR{X_k>W_k}=1\qquad\quad\frac{a\mu}{\lambda}>1
\end{equation*}
On the other hand for $^{a\mu}/_\lambda\le1$ the probability
$\PR{X_k>W_k}$ can still be explicitly calculated by taking into
account the laws specified in the Section\myref{pcpp}, and in this
case it is possible to show that
\begin{equation*}
    \PR{X_k>W_k}=\frac{(1-a)\mu}{\lambda+(1-a)\mu}\qquad\quad\frac{a\mu}{\lambda}\le1
\end{equation*}
a value ranging from $0$ to $1$ according to the different possible
choices of the parameters $a,\lambda$ and $\mu$

As for the relative timings $T_n,S_m$ of the shocks along the point
processes themselves, an explicit calculation of $\PR{T_n\le S_m}$
is certainly possible, but its results would turn out to be rather
cumbersome because it would involve two or three convolutions of
(positive and negative) Erlang laws with different parameters. We
will then confine ourselves here to produce just the
cross-correlations between $T_n, S_m$: since it is easy to check
that
\begin{equation*}
    \cov{X_k}{W_\ell}=\frac{a}{\lambda\mu}\,\delta_{k\ell}
\end{equation*}
it is also apparent that
\begin{equation*}
    \cov{T_n}{S_m}=\sum_{k=1}^{n}\sum_{\ell=1}^{m}\cov{X_k}{W_\ell}=\frac{a}{\lambda\mu}\sum_{k=1}^{n}\sum_{\ell=1}^{m}\delta_{k\ell}=\frac{a}{\lambda\mu}m\wedge
    n
\end{equation*}
and hence the cross-correlation coefficient of $T_n, S_m$ will
simply be
\begin{equation*}
    r_{nm}=a\,\frac{n\wedge m}{\sqrt{nm}}=\left\{
                                             \begin{array}{ll}
                                               a\sqrt{\,^n/_m} &\quad \hbox{for $\;n\le m$} \\
                                               a\sqrt{\,^m/_n} &\quad \hbox{for $\;n\ge m$}
                                             \end{array}
                                           \right.
\end{equation*}
\begin{figure}
\begin{center}
\includegraphics*[width=10cm]{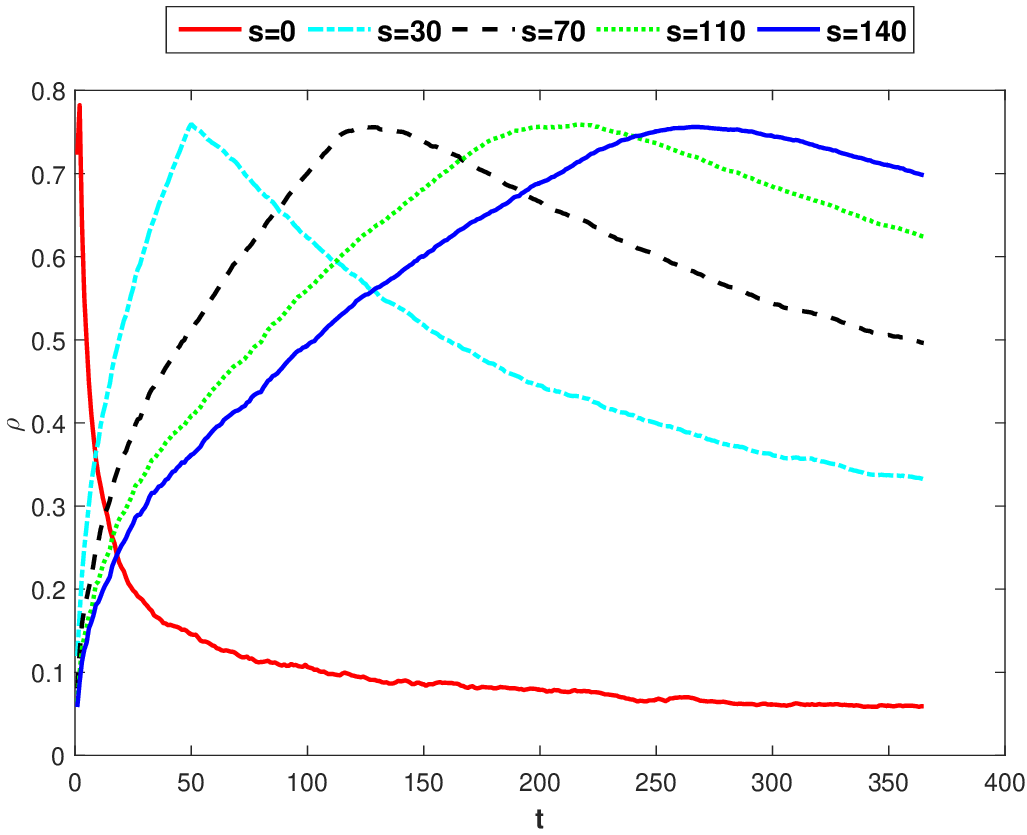}
\caption{Cross-correlation $\rho(s,t)$ of the two Poisson processes
$M(s)$ and $N(t)$ estimated by Monte Carlo simulations. Here
$a=0.8,\;\lambda=20$ and $\mu=40$, namely
$^{a\mu}/_\lambda=1.6>1$}\label{crosscorr08}
\end{center}
\end{figure}

Finally, even closed formulas for the cross-correlation coefficient
$\rho(s,t)$ between the Poisson processes $M(s)$ and $N(t)$ are
still derivable on the ground of our previous results about the
joint distributions, but it would be too long to thoroughly
elaborate them here. As an alternative we have chosen to show the
plots of their estimates based on a sample of $10^5$ Monte Carlo
simulations of their trajectories as shown in the
Figures\myref{crosscorr05} and\myref{crosscorr08}. There the
behavior is displayed of $\rho(s,t)$ as a function of $t$ for
different, fixed values of $s$. More precisely, in the
Figure\myref{crosscorr05} we have taken $a=0.5$ and $\lambda=\mu=20$
as the values for the relevant parameters of our coupled processes
(then we have $^{a\mu}/_\lambda<1$), while in the
Figure\myref{crosscorr08} the same parameters are $a=0.8,\lambda=20$
and $\mu=40$ (and then $^{a\mu}/_\lambda>1$). It is apparent from
these pictures that the behavior of $\rho(s,t)$ is comparable to
that of the self-correlation of a simple Poisson process, but for
the fact that the cumulative effect of the correlate renewals
produces a smoothing of the shape around the maximum values near
$t=s$. At first sight this could look as a little difference, but in
the domain, for instance, of the financial applications even small
deviations among the models could produce huge differences in gains
and losses

\section{Conclusions and further inquiries}\label{negative}

It is apparent that, within the model discussed in the
Section\myref{sdlaws}, from the self-decomposability alone we can
only get pairs of \rv's $X,Y$ with \emph{positive correlations}
$0<r_{XY}<1$ steered by the value of a parameter $a$. It would be
interesting however to widen the scope of our models in order to
achieve also Poisson processes whose correlation can span over all
its possible values (both positive and negative) by changing the
value of some numerical parameter. In this respect it is important
to remark -- as pointed out in the Section\myref{frechet} -- that
while two \rv's $X$ and $Y$ which are, for instance, marginally
exponentials can also be \emph{totally correlated} ($r_{XY}=1$),
they can not instead be \emph{totally anti-correlated} ($r_{XY}=-1$)
because this would imply some \emph{linear} dependence with a
\emph{negative} proportionality coefficient, and that would be at
odds with the fact that both our \rv's take arbitrary large, but
only positive values. Hence two exponential \rv's $X$ and $Y$ can
always have a negative correlation, but only up to a minimal value
which in any case must be larger than $-1$. We also showed in the
Section\myref{frechet} that this minimum is reached when between $X$
and $Y$ there is a peculiar kind of mutual functional,
\emph{decreasing} dependence, albeit clearly \emph{not a linear}
one. A model to produce pairs $X,Y$ of \rv's which are marginally
exponentials, and which -- following the value of a numerical
parameter $a$ -- show all the possible correlation values will be
discussed in a subsequent paper

Our results in any case show that the self-decomposability, joined
with the subordination techniques, can be a promising tool to study
dependency beyond the Gaussian-{\ito} world. We have shown indeed
how to obtain dependent exponential (gamma) \rv's that can be used
to create and simulate dependent Poisson processes without resorting
to definitely coincident jumps, but the path is now open to produce
more general dependent gamma (Erlang at first) \rv's to simulate
dependent variance gamma processes. A further extension could then
be to study the self-decomposability of density functions that have
a strictly proper rational characteristic function (Exponential
Polynomial Trigonometric, \emph{EPT} laws) in order to construct
2-dimensional correlated \emph{EPT} \rv's (see for
instance\mycite{SextonHanzon12}). Finally it would be expedient to
explore the Markov properties of the 2-component Poisson processes
$(M,(t),N(t))$ with dependent marginals that we have introduced in
this paper and the Master equations ruling them: this too will be
the subject of future inquiries

% ------------------------------------------------------------------------

\subsection*{Acknowledgment}
The authors would like to thank prof.\ R.\ Nelsen and prof.\ S.\
Stramaglia for invaluable discussions and suggestions

\begin{appendix}

\section{Notations}\label{notations}

All along this paper we will adopt the following notations: for a
Poisson law $\poiss(\alpha)$ we will introduce the symbols
\begin{equation*}
    \pi_n(\alpha)=\frac{\alpha^n}{n!}\,e^{-\alpha}\qquad\quad\Pi_n(\alpha)=\sum_{k=n}^{\infty}\pi_k(\alpha)\quad\qquad\alpha>0\qquad
    n=0,1,2\ldots
\end{equation*}
and for a binomial law $\bin(n,1-a)$ the notation
\begin{equation*}
    \beta_k(n)=\left\{
                 \begin{array}{cll}
                   1 &\quad n=0, &k=0 \\
                   \binom{n}{k}a^{n-k}(1-a)^k &\quad
                   n=1,2,\ldots, &k=0,\ldots,n
                 \end{array}
               \right.
    \qquad 0\le a\le 1
\end{equation*}
It will be understood moreover that
\begin{equation*}
    \pi_n(0^+)=\delta_{n,0}
\end{equation*}
We will also use bot the Heaviside function $\heavi$, and the
Heaviside symbol $\Theta$
\begin{equation*}
    \heavi(x)=\left\{
                   \begin{array}{ll}
                     1 & \quad x\ge0 \\
                     0 & \quad x<0
                   \end{array}
                 \right.
                 \qquad\qquad
    \Theta_j=\left\{
               \begin{array}{ll}
                 1 & \quad j\ge1 \\
                 0 & \quad j\le0
               \end{array}
             \right.\qquad j=0,\pm1,\pm2,\ldots
\end{equation*}
The \pdf\ and \chf\ of a standard Erlang law $\erl_n(1)$ moreover
will be denoted as
\begin{equation*}
    f_n(x)=\left\{
             \begin{array}{c}
               \delta(x) \\
               \frac{x^{n-1}}{(n-1)!}e^{-x}\heavi(x)
             \end{array}
           \right.
           \qquad
    \varphi_n(u)=\left\{
             \begin{array}{cl}
               1 & \qquad n=0 \\
               \left(\frac{1}{1-iu}\right)^n & \qquad n=1,2,\ldots
             \end{array}
           \right.
\end{equation*}
where it is understood for the Dirac delta $\delta(x)$ that for
every $b>0$
\begin{equation*}
    \int_0^b\delta(x)\,dx=1\qquad\qquad\lim_{z\downarrow0^+}\int_z^b\delta(x)\,dx=0
\end{equation*}
Remark also that apparently
\begin{equation*}
    f_n(x)=\pi_{n-1}(x)\vartheta(x)\qquad n=1,2,\ldots
\end{equation*}
We will finally define for later convenience the functions
\begin{equation*}
    h_n(x)=\sum_{k=0}^n\beta_k(n)f_k(x)\qquad n=0,1,2,\ldots
\end{equation*}
which are the \pdf's of the mixtures of Erlang laws $\erl_k(1)$ with
binomial $\bin(n,1-a)$ weights for their indices $k$

\section{A proof of Proposition\myref{prop}}\label{proofprop}

To evaluate $\PR{S_m\le \rho,\,S_n\le \tau}$ we first remark that
\begin{equation*}
  \PR{S_m\le \rho,\,S_n\le \tau} = \PR{M(\rho)\ge m,\,M(\tau)\ge n}
\end{equation*}
and then that, being a Poisson process, $M(t)$ is non-decreasing: as
a consequence \begin{eqnarray*}
    m\le n\;\;\hbox{and}\;\; \tau\le\rho &\Longrightarrow& M(\tau)\le M(\rho)\;\;\;\hbox{hence}\;\;\{M(\tau)\ge n\}\subseteq\{M(\rho)\ge m\}\\
    n\le m\;\;\hbox{and}\;\;\rho\le\tau &\Longrightarrow& M(\rho)\le M(\tau)\;\;\;\hbox{hence}\;\;\{M(\rho)\ge m\}\subseteq\{M(\tau)\ge n\}
\end{eqnarray*}
In the case $m\le n$ we then have for $ \tau\le\rho$
\begin{equation*}
   \PR{M(\rho)\ge m,\,M(\tau)\ge n}=\PR{M(\tau)\ge n}
\end{equation*}
while for $\rho\le\tau$ from the general properties of a Poisson
process we get
\begin{eqnarray*}
  \lefteqn{\PR{M(\rho)\ge m,\,M(\tau)\ge n}}\qquad\quad\\
   &=& \sum_{k=m}^\infty\PR{M(\rho)\ge m,\,M(\tau)\ge n\,|\,M(\rho)=k}\PR{M(\rho)=k} \\
   &=&\sum_{k=m}^n\PR{M(\tau)\ge n\,|\,M(\rho)=k}\PR{M(\rho)=k}+\PR{M(\rho)>n} \\
   &=&\sum_{k=m}^n\PR{M( \tau-\rho)\ge n-k}\PR{M(\rho)=k}+\PR{M(\rho)>n}
\end{eqnarray*}
In the same vein, when $n\le m$ we have for $\rho\le\tau$
\begin{equation*}
   \PR{M(\rho)\ge m,\,M(\tau)\ge n}=\PR{M(\rho)\ge n}
\end{equation*}
while for $ \tau\le\rho$ we get
\begin{eqnarray*}
  \lefteqn{\PR{M(\rho)\ge m,\,M(\tau)\ge n}}\qquad\quad\\
   &=& \sum_{k=m}^\infty\PR{M(\rho)\ge m,\,M(\tau)\ge n\,|\,M(\tau)=k}\PR{M(\tau)=k} \\
   &=&\sum_{k=n}^m\PR{M(\rho)\ge m\,|\,M(\tau)=k}\PR{M(\tau)=k}+\PR{M(\tau)>m} \\
   &=&\sum_{k=n}^m\PR{M(\rho-\tau)\ge m-k}\PR{M(\tau)=k}+\PR{M(\tau)>m}
\end{eqnarray*}
Remark that for $m=n$ both the cases lead to the same result, namely
\begin{equation*}
    \PR{M(\rho)\ge n,\,M(\tau)\ge n}=\left\{
                                    \begin{array}{ll}
                                      \PR{M(\tau)\ge n} & \quad\hbox{when} \; \tau\le\rho \\
                                      \PR{M(\rho)\ge n} & \quad\hbox{when} \;\rho\le\tau
                                    \end{array}
                                  \right.
\end{equation*}
that can also be conveniently summarized as
\begin{equation*}
    \PR{M(\rho)\ge n,\,M(\tau)\ge n}=\PR{M(\rho\wedge\tau)\ge n}
\end{equation*}
On the other hand for $m<n$ we have
\begin{eqnarray*}
  \hbox{for}\;\; \tau\le\rho &&\PR{M(\tau)\ge n}  \\
  \hbox{for}\;\;\rho\le\tau &&\PR{M(\rho)\ge n}+\sum_{k=m}^{n-1}\PR{M( \tau-\rho)\ge n-k}\PR{M(\rho)=k}
\end{eqnarray*}
that can also be put in the form
\begin{equation*}
  \PR{M(\rho\wedge\tau)\ge n}+\heavi( \tau-\rho)\sum_{k=m}^{n-1}\PR{M( \tau-\rho)\ge n-k}\PR{M(\rho)=k}
\end{equation*}
while for $m>n$ it is
\begin{eqnarray*}
  \hbox{for}\;\; \tau\le\rho && \PR{M(\tau)\ge m}+\sum_{k=n}^{m-1}\PR{M(\rho-\tau)\ge m-k}\PR{M(\tau)=k} \\
  \hbox{for}\;\;\rho\le\tau &&\PR{M(\rho)\ge m}
\end{eqnarray*}
namely
\begin{equation*}
  \PR{M(\rho\wedge\tau)\ge m}+\heavi(\rho-\tau)\sum_{k=n}^{m-1}\PR{M(\rho-\tau)\ge m-k}\PR{M(\tau)=k}
\end{equation*}
In both cases the first terms can expressed as
$\PR{M(\rho\wedge\tau)\ge m\vee n}$, and in this form they also
coincide with the previous result for $m=n$. On the other hand the
extra term with the sum (which is absent for $m=n$) must be taken in
consideration either when we have both $m<n$ and $\rho\le\tau$, or
when it is $m>n$ and $ \tau\le\rho$. All these provisions can then
be comprehensively taken into account in the formula
\begin{eqnarray*}
  \lefteqn{\PR{S_m\le \rho,\,S_n\le \tau}}\\
   &=&\!\! \PR{M(\rho\wedge \tau)\ge m\vee n}+\left[\Theta_{n-m}\heavi(\tau-\rho)+\Theta_{m-n}\heavi(\rho-\tau)\right]\cdot \\
   &&\qquad\qquad\sum_{k=m\wedge
   n}^{(m\vee n)-1}\!\!\PR{M(|\rho-\tau|)\ge(m\vee n)-k}\PR{M(\rho\wedge \tau)=k}
\end{eqnarray*}
which finally takes the form of Proposition\myref{prop} by using the
notations adopted in the Appendix\myref{notations} for the Poisson
distributions and the Heaviside symbols

\section{Proof details for Proposition\myref{case1}}\label{proof1}

\subsection{The case $a\mu s>\lambda t$, namely $z<0$}

We begin with the case $a\mu s>\lambda t$, namely $z<0$, by
recalling also that for the integration variable in\refeq{qmn} it is
$0\le w\le y$. As a consequence, when the Proposition\myref{prop} in
used in\refeq{qmn}, we always have
\begin{equation}\label{ineq}
    0\le \tau=\frac{y-w}{\mu}\le\frac{y-z}{\mu}=s=\rho
\end{equation}
On the other hand, since the conditions of the Lemma\myref{lemma}
are met, we can also restrict ourselves to evaluate $p_{m,n}(s,t)$
for $0\le n\le m$

\subsubsection{$\bm{m=n\ge0}$}

In this case from\refeq{qmn},\refeq{ineq} and from the
Proposition\myref{prop} we find
\begin{eqnarray*}
  q_{n,n}(y,z)&=&a\int_0^y\!\!dw\,h_n(aw)\PR{M\left(\frac{y-w}{\mu}\right)\ge n} \\
  q_{n+1,n}(y,z) &=&a\int_0^y\!\!dw\,h_n(aw)\bigg[\PR{M\left(\frac{y-w}{\mu}\right)\ge
  n+1}\\
     &&\qquad\qquad\qquad+\,\PR{M\left(\frac{w-z}{\mu}\right)\ge1}\PR{M\left(\frac{y-w}{\mu}\right)=n}\bigg] \\
  q_{n,n+1}(y,z) &=&a\int_0^y\!\!dw\,h_{n+1}(aw)\PR{M\left(\frac{y-w}{\mu}\right)\ge n+1} \\
  q_{n+1,n+1}(y,z) &=&a\int_0^y\!\!dw\,h_{n+1}(aw)\PR{M\left(\frac{y-w}{\mu}\right)\ge n+1}
\end{eqnarray*}
and then from\refeq{pmn}
\begin{eqnarray}
  p_{n,n}(y,z)\!\! &=&\!\!a\int_0^y\!\!dw\,h_n(aw)\PR{M\left(\frac{w-z}{\mu}\right)=0}\PR{M\left(\frac{y-w}{\mu}\right)=n}\nonumber \\
   &=&\!\!a\int_0^y\!\!dw\,h_n(aw)\,\pi_0(w-z)\,\pi_n(y-w)\label{z<0m=n}
\end{eqnarray}

\vfill\eject

\subsubsection{$\bm{m>n\ge0}$}

In order to make use of Proposition\myref{prop} we remark that now
\begin{equation*}
    m\ge n+1\qquad\quad m+1>n\qquad\quad m+1>n+1
\end{equation*}
so that from\refeq{qmn} and\refeq{ineq}
\begin{eqnarray*}
  q_{m,n}(y,z)\!\!\! &=& \!\!\!a\int_0^y\!\!dw\,h_n(aw)\Bigg[\PR{M\left(\frac{y-w}{\mu}\right)\ge m}\Bigg.\\
                &&  \quad       \left.+\sum_{k=n}^{m-1}\PR{M\left(\frac{w-z}{\mu}\right)\ge m-k}\PR{M\left(\frac{y-w}{\mu}\right)=k}\right] \\
  q_{m+1,n}(y,z)\!\!\! &=& \!\!\!a\int_0^y\!\!dw\,h_n(aw)\Bigg[\PR{M\left(\frac{y-w}{\mu}\right)\ge m+1}\Bigg.\\
                &&  \quad       \left.+\sum_{k=n}^{m}\PR{M\left(\frac{w-z}{\mu}\right)\ge m+1-k}\PR{M\left(\frac{y-w}{\mu}\right)=k}\right] \\
  q_{m,n+1}(y,z) \!\!\!&=&\!\!\!a\int_0^y\!\!dw\,h_{n+1}(aw)\Bigg[\PR{M\left(\frac{y-w}{\mu}\right)\ge m}\Bigg.\\
                &&     \left.+\,\Theta_{m-n-1}\!\!\sum_{k=n+1}^{m-1}\!\!\PR{M\left(\frac{w-z}{\mu}\right)\ge m-k}\PR{M\left(\frac{y-w}{\mu}\right)=k}\right] \\
  q_{m+1,n+1}(y,z)\!\!\! &=& \!\!\!a\int_0^y\!\!dw\,h_{n+1}(aw)\Bigg[\PR{M\left(\frac{y-w}{\mu}\right)\ge m+1}\Bigg.\\
                &&  \quad       \left.+\sum_{k=n+1}^{m}\PR{M\left(\frac{w-z}{\mu}\right)\ge m+1-k}\PR{M\left(\frac{y-w}{\mu}\right)=k}\right]
\end{eqnarray*}
and hence from\refeq{pmn} (in a slightly simplified notation)
\begin{eqnarray*}
  p_{m,n}\!\! &=&\!\!a\int_0^y\!\!dw\,h_n\Bigg[\PR{M\left(\frac{y-w}{\mu}\right)=m}\\
                && \qquad\qquad\quad-\PR{M\left(\frac{w-z}{\mu}\right)\ge 1}\PR{M\left(\frac{y-w}{\mu}\right)=m}\\
                &&  \qquad\qquad\qquad       \left.+\sum_{k=n}^{m-1}\PR{M\left(\frac{w-z}{\mu}\right)=m-k}\PR{M\left(\frac{y-w}{\mu}\right)=k}\right] \\
                && -a\int_0^y\!\!dw\,h_{n+1}\Bigg[\PR{M\left(\frac{y-w}{\mu}\right)=m}\\
                && \qquad-\PR{M\left(\frac{w-z}{\mu}\right)\ge 1}\PR{M\left(\frac{y-w}{\mu}\right)=m}\\
                &&  \qquad\qquad
                \left.+\,\Theta_{m-n-1}\sum_{k=n+1}^{m-1}\PR{M\left(\frac{w-z}{\mu}\right)=m-k}\PR{M\left(\frac{y-w}{\mu}\right)=k}\right]
                \\
            &=& \!\!a\int_0^y\!\!dw\Bigg[(h_n-h_{n+1})\PR{M\left(\frac{w-z}{\mu}\right)=0}\PR{M\left(\frac{y-w}{\mu}\right)=m}
            \\
            &&\qquad\qquad\qquad +h_n\PR{M\left(\frac{w-z}{\mu}\right)=m-n}\PR{M\left(\frac{y-w}{\mu}\right)=n} \\
            &&\left.+(h_n-h_{n+1})\Theta_{m-n-1}\!\!\sum_{k=n+1}^{m-1}\!\!\PR{M\left(\frac{w-z}{\mu}\right)=m-k}\PR{M\left(\frac{y-w}{\mu}\right)=k}\right]\\
            &=&\!\!a\int_0^y\!\!dw\Bigg[h_{n+1}\PR{M\left(\frac{w-z}{\mu}\right)=m-n}\PR{M\left(\frac{y-w}{\mu}\right)=n}\\
            &&\qquad\left.+(h_n-h_{n+1})\sum_{k=n}^{m}\PR{M\left(\frac{w-z}{\mu}\right)=m-k}\PR{M\left(\frac{y-w}{\mu}\right)=k}\right]
\end{eqnarray*}
so that we finally have in full notation
\begin{eqnarray}
  p_{m,n}(y,z)&=&a\int_0^y\!\!dw\Bigg\{h_{n+1}(aw)\pi_{m-n}(w-z)\pi_n(y-w)\label{z<0m>n} \\
    &&\qquad\qquad\qquad+\big[h_n(aw)-h_{n+1}(aw)\big]\sum_{k=n}^{m}\pi_{m-k}(w-z)\pi_k(y-w)\Bigg\}\nonumber
\end{eqnarray}
Remark that this formula correctly encompasses also the case
$m=n\ge0$. As stated in the Proposition\myref{case1}, we can also
concisely write the overall result in the form
\begin{equation}\label{pmn1}
    p_{m,n}(y,z)=\left\{
                   \begin{array}{ll}
                     Q_{n,n}(y,z) & \qquad m=n\ge0 \\
                     Q_{m,n}(y,z)-Q_{m,n+1}(y,z) & \qquad m>n\ge0
                   \end{array}
                 \right.
\end{equation}
where for $m\ge n\ge0$ and $z<0$ we define
\begin{equation*}
  Q_{m,n}(y,z)= a\int_0^ydw\, h_n(aw)\sum_{k=n}^{m}\pi_{m-k}(w-z)\pi_k(y-w)
\end{equation*}

\subsubsection{A normalization check}\label{normchk} We provide here
a quick check of the correct normalization of the joint
probabilities $p_{m,n}$ calculated in the previous section for $a\mu
s>\lambda t$, namely $z<0<y$. First remark that with $x=aw$
\begin{eqnarray*}
  Q_{m,0}(y,z) &=& a\int_0^ydw\, h_0(aw)\sum_{k=0}^{m}\pi_{m-k}(w-z)\pi_k(y-w) \\
   &=&\int_0^{ay}dx\,
   \delta(x)\sum_{k=0}^{m}\pi_{m-k}\left(\frac{x}{a}-z\right)\pi_k\left(y-\frac{x}{a}\right)=\sum_{k=0}^{m}\pi_{m-k}(-z)\pi_k(y)\\
   &=&e^{-(y-z)}\sum_{k=0}^{m}\frac{y^k(-z)^{m-k}}{k!(m-k)!}=\frac{e^{-(y-z)}}{m!}\sum_{k=0}^{m}\binom{m}{k}y^k(-z)^{m-k}\\
   &=&\frac{e^{-(y-z)}}{m!}(y-z)^m=\pi_m(y-z)
\end{eqnarray*}
Then, since here only the $p_{m,n}$ with $m\ge n\ge0$ do not vanish,
we have (by neglecting the arguments $y,z$)
\begin{eqnarray*}
    \sum_{n,m}p_{m,n}&=&\sum_{m=0}^\infty\sum_{n=0}^mp_{m,n}=Q_{0,0}+\sum_{m=1}^\infty\left[Q_{m,m}+\sum_{n=0}^{m-1}\left(Q_{m,n}-Q_{m,n+1}\right)\right]\\
    &=&Q_{0,0}+\sum_{m=1}^\infty\left(\sum_{n=0}^mQ_{m,n}-\sum_{n=0}^{m-1}Q_{m,n+1}\right)\\
    &=&Q_{0,0}+\sum_{m=1}^\infty\left(\sum_{n=0}^mQ_{m,n}-\sum_{n=1}^mQ_{m,n}\right)=Q_{0,0}+\sum_{m=1}^\infty
    Q_{m,0}\\
    &=&\sum_{m=0}^\infty Q_{m,0}=\sum_{m=0}^\infty\pi_m=1
\end{eqnarray*}
which confirms the normalization

\subsection{The case $a\mu s<\lambda t$, namely $y>z>0$}

In this case, since $0<z<y$ and $0<w<y$, when we use
Proposition\myref{prop} in\refeq{qmn} instead of\refeq{ineq} we find
only
\begin{equation}\label{ineq2}
   0\le \tau=\frac{y-w}{\mu}\qquad\quad0\le \rho=s=\frac{y-z}{\mu}
\end{equation}
so that $\rho$ and $\tau$ can now happen to be in an order
whatsoever. As a consequence Lemma\myref{lemma} does not hold, and
we must consider all the possible choices of $m,n$

\subsubsection{$\bm{n>m\ge0}$}

In order to make use of Proposition\myref{prop} we remark that now
\begin{equation*}
    m< n+1\qquad\quad m+1\le n\qquad\quad m+1<n+1
\end{equation*}
so that from\refeq{qmn} and\refeq{ineq2}
\begin{eqnarray*}
  q_{m,n}\!\! &=& \!\!a\int_0^y\!\!dw\,h_n(aw)\PR{M\left(\frac{y-z}{\mu}\,\wedge\,\frac{y-w}{\mu}\right)\ge n}\\
                &&   \qquad  +a\int_0^z\!\!dw\,h_n(aw)\\
                && \qquad\qquad\sum_{k=m}^{n-1}\PR{M\left(\frac{z-w}{\mu}\right)\ge n-k}\PR{M\left(\frac{y-z}{\mu}\right)=k} \\
  q_{m+1,n}\!\! &=& \!\!a\int_0^y\!\!dw\,h_n(aw)\PR{M\left(\frac{y-z}{\mu}\,\wedge\,\frac{y-w}{\mu}\right)\ge n}\\
                &&  \qquad    +\,a\,\Theta_{n-m-1}\int_0^z\!\!dw\,h_n(aw)\\
                && \qquad\qquad\sum_{k=m+1}^{n-1}\PR{M\left(\frac{z-w}{\mu}\right)\ge n-k}\PR{M\left(\frac{y-z}{\mu}\right)=k} \\
  q_{m,n+1} \!\!&=&\!\!  \!\!a\int_0^y\!\!dw\,h_{n+1}(aw)\PR{M\left(\frac{y-z}{\mu}\,\wedge\,\frac{y-w}{\mu}\right)\ge n+1}\\
                &&  \qquad      +\,a\int_0^z\!\!dw\,h_{n+1}(aw)\\
                && \qquad\qquad\sum_{k=m}^{n}\PR{M\left(\frac{z-w}{\mu}\right)\ge n+1-k}\PR{M\left(\frac{y-z}{\mu}\right)=k} \\
  q_{m+1,n+1}\!\! &=&  \!\!a\int_0^y\!\!dw\,h_{n+1}(aw)\PR{M\left(\frac{y-z}{\mu}\,\wedge\,\frac{y-w}{\mu}\right)\ge n+1}\\
                &&  \qquad      +\,a\int_0^z\!\!dw\,h_{n+1}(aw)\\
                && \qquad\qquad\sum_{k=m+1}^{n}\PR{M\left(\frac{z-w}{\mu}\right)\ge n+1-k}\PR{M\left(\frac{y-z}{\mu}\right)=k}
\end{eqnarray*}
and hence from\refeq{pmn}
\begin{eqnarray*}
  p_{m,n} \!\!&=& a\int_0^z\!\!dw\,h_n(aw)\PR{M\left(\frac{z-w}{\mu}\right)\ge n-m}\PR{M\left(\frac{y-z}{\mu}\right)=m} \\
   &&\!\!-a\int_0^z\!\!dw\,h_{n+1}(aw)\PR{M\left(\frac{z-w}{\mu}\right)\ge
   n+1-m}\PR{M\left(\frac{y-z}{\mu}\right)=m} \\
   &=&\!\!\!\PR{M\left(\frac{y-z}{\mu}\right)=m}a\!\!\int_0^z\!\!dw\,\Bigg[h_n(aw)\left(1-\PR{M\left(\frac{z-w}{\mu}\right)<n-m}\right)\\
   && \qquad\qquad\qquad\qquad -h_{n+1}(aw)\left(1-\PR{M\left(\frac{z-w}{\mu}\right)\le
   n-m}\right)\Bigg]\\
   &=&\PR{M\left(\frac{y-z}{\mu}\right)=m}a\int_0^z\!\!dw\,\Bigg\{h_n(aw)\PR{M\left(\frac{z-w}{\mu}\right)=n-m}\\
   && \qquad\qquad\qquad +[h_n(aw)-h_{n+1}(aw)]\left(1-\PR{M\left(\frac{z-w}{\mu}\right)\le
   n-m}\right)\Bigg\}
\end{eqnarray*}
which, by plugging in the explicit Poisson probabilities, can be
written as
\begin{eqnarray}
  p_{m,n}(y,z) &=& \pi_m(y-z)\,a\int_0^z\!\!dw\,[h_n(aw)-h_{n+1}(aw)] \nonumber\\
   &&\qquad+\pi_m(y-z)\,a\int_0^z\!\!dw\Big\{h_n(aw)\pi_{n-m}(z-w)\label{z>0n>m}\\
   &&\qquad\qquad\qquad\qquad\qquad-[h_n(aw)-h_{n+1}(aw)]\sum_{k=0}^{n-m}\pi_k(z-w)\Big\}\nonumber
   \\
   &=&A_{m,n}(y,z)-A_{m,n+1}(y,z)+B_{m,n}(y,z)-B_{m,n-1}(y,z)\nonumber
\end{eqnarray}
where, with the notations\refeq{A1} and\refeq{B1} adopted in the
Proposition\myref{case1}, we have defined
\begin{eqnarray*}
  A_{m,n}(y,z)&=& \pi_m(y-z)\;a\int_0^z\!\!dw\,h_n(aw)\\
  B_{m,n}(y,z)&=& \pi_m(y-z)\;a\int_0^z\!\!dw\,h_{n+1}(aw)\sum_{k=0}^{n-m}\pi_k(z-w)
\end{eqnarray*}

\subsubsection{$\bm{n=m\ge0}$}

In this case from\refeq{qmn},\refeq{ineq2} and from the
Proposition\myref{prop} we find
\begin{eqnarray*}
  q_{n,n}\!\! &=& \!\!a\int_0^y\!\!dw\,h_n(aw)\PR{M\left(\frac{y-z}{\mu}\wedge\frac{y-w}{\mu}\right)\ge n} \\
  q_{n+1,n}\!\! &=& \!\!a\int_0^y\!\!dw\,h_n(aw)\PR{M\left(\frac{y-z}{\mu}\wedge\frac{y-w}{\mu}\right)\ge n+1}\\
  &&\quad+a\int_z^y\!\!dw\,h_n(aw)\PR{M\left(\frac{w-z}{\mu}\right)\ge1}\PR{M\left(\frac{y-w}{\mu}\right)=n} \\
  q_{n,n+1} \!\!&=&\!\! a\int_0^y\!\!dw\,h_{n+1}(aw)\PR{M\left(\frac{y-z}{\mu}\wedge\frac{y-w}{\mu}\right)\ge n+1} \\
  &&\quad+a\int_0^z\!\!dw\,h_{n+1}(aw)\PR{M\left(\frac{w-z}{\mu}\right)\ge1}\PR{M\left(\frac{y-w}{\mu}\right)=n} \\
  q_{n+1,n+1}\!\! &=& \!\! a\int_0^y\!\!dw\,h_{n+1}(aw)\PR{M\left(\frac{y-z}{\mu}\wedge\frac{y-w}{\mu}\right)\ge n+1}
\end{eqnarray*}
and then from\refeq{pmn}
\begin{eqnarray*}
  p_{m,n} \!\! &=& \!\!a\int_0^y\!\!dw\,h_n(aw)\PR{M\left(\frac{y-z}{\mu}\wedge\frac{y-w}{\mu}\right)=n}\\
   &&\;-a\int_z^y\!\!dw\,h_n(aw)\left(1-\PR{M\left(\frac{w-z}{\mu}\right)=0}\right)\PR{M\left(\frac{y-w}{\mu}\right)=n}\\
   &&-a\!\int_0^z\!\!dw\,h_{n+1}(aw)\!\left(1-\PR{M\left(\frac{w-z}{\mu}\right)=0}\right)\!\PR{M\left(\frac{y-w}{\mu}\right)=n}\\
   &=& \!\!a\int_0^z\!\!dw\,[h_n(aw)-h_{n+1}(aw)]\,\PR{M\left(\frac{y-z}{\mu}\right)=n}\\
   &&\qquad-a\int_z^y\!\!dw\,h_n(aw)\PR{M\left(\frac{w-z}{\mu}\right)=0}\PR{M\left(\frac{y-w}{\mu}\right)=n}\\
   &&\qquad-a\int_0^z\!\!dw\,h_{n+1}(aw)\PR{M\left(\frac{w-z}{\mu}\right)=0}\PR{M\left(\frac{y-w}{\mu}\right)=n}
\end{eqnarray*}
which with the Poisson probabilities become
\begin{eqnarray}
  p_{n,n}(y,z) &=& \pi_n(y-z)\,a\int_0^z\!\!dw\,[h_n(aw)-h_{n+1}(aw)]\nonumber\\
   &&\qquad\qquad+\pi_n(y-z)\,a\int_0^z\!\!dw\, h_{n+1}(aw)\pi_0(z-w)\label{z>0m=n} \\
   &&\qquad\qquad\qquad\qquad\qquad+a\int_z^y\!\!dw \,h_n(aw)\pi_0(w-z)\pi_n(y-w) \nonumber\\
   &=&A_{n,n}(y,z)-A_{n,n+1}(y,z)+B_{n,n}(y,z)+C_{n,n}(y,z)\nonumber
\end{eqnarray}
where, as in the Proposition\myref{case1}, along with\refeq{A1}
and\refeq{B1} we also introduced\refeq{C1}
\begin{equation*}
    C_{m,n}(y,z)= a\int_z^y\!\!dw\,h_n(aw)\sum_{k=n}^m\pi_{m-k}(w-z)\pi_k(y-w)
\end{equation*}

\subsubsection{$\bm{m>n\ge0}$}

In this case it is
\begin{equation*}
    m\ge n+1\qquad\quad m+1> n\qquad\quad m+1>n+1
\end{equation*}
so that from\refeq{qmn} and\refeq{ineq2}
\begin{eqnarray*}
  q_{m,n}\!\! &=& \!\!a\int_0^y\!\!dw\,h_n(aw)\PR{M\left(\frac{y-z}{\mu}\wedge\frac{y-w}{\mu}\right)\ge m}\\
                && \quad +\,a\int_z^y\!\!dw\,h_n(aw)\\
                && \qquad\quad\sum_{k=n}^{m-1}\PR{M\left(\frac{w-z}{\mu}\right)\ge m-k}\PR{M\left(\frac{y-w}{\mu}\right)=k} \\
  q_{m+1,n}\!\! &=& \!\!a\int_0^y\!\!dw\,h_n(aw)\PR{M\left(\frac{y-z}{\mu}\wedge\frac{y-w}{\mu}\right)\ge m+1}\\
                &&  \quad      +\,a\int_z^y\!\!dw\,h_n(aw)\\
                &&\qquad\quad\sum_{k=n}^{m}\PR{M\left(\frac{w-z}{\mu}\right)\ge m+1-k}\PR{M\left(\frac{y-w}{\mu}\right)=k} \\
  q_{m,n+1} \!\!&=&\!\!  \!\!a\int_0^y\!\!dw\,h_{n+1}(aw)\PR{M\left(\frac{y-z}{\mu}\wedge\frac{y-w}{\mu}\right)\ge m}\\
                &&  \quad+\,\Theta_{m-n-1}a\int_z^y\!\!dw\,h_{n+1}(aw)\\
                &&\qquad\quad\sum_{k=n+1}^{m-1}\PR{M\left(\frac{w-z}{\mu}\right)\ge m-k}\PR{M\left(\frac{y-w}{\mu}\right)=k} \\
  q_{m+1,n+1}\!\! &=&  \!\!a\int_0^y\!\!dw\,h_{n+1}(aw)\PR{M\left(\frac{y-z}{\mu}\wedge\frac{y-w}{\mu}\right)\ge m+1}\\
                &&  \quad      +\,a\int_z^y\!\!dw\,h_{n+1}(aw)\\
                &&\qquad\quad\sum_{k=n+1}^{m}\PR{M\left(\frac{w-z}{\mu}\right)\ge m+1-k}\PR{M\left(\frac{y-w}{\mu}\right)=k}
\end{eqnarray*}
and hence from\refeq{pmn}
\begin{eqnarray*}
   p_{m,n} \!\! &=&\!\!a\int_0^y\!\!dw\,[h_n(aw)-h_{n+1}(aw)]\,\PR{M\left(\frac{y-z}{\mu}\wedge\frac{y-w}{\mu}\right)=m}  \\
   &&+a\int_z^y\!\!dw\,h_n(aw)\left[\sum_{k=n}^{m-1}\PR{M\left(\frac{w-z}{\mu}\right)=m-k}\PR{M\left(\frac{y-w}{\mu}\right)=k}\right.\\
   &&\qquad\qquad\qquad\qquad-\PR{M\left(\frac{w-z}{\mu}\right)\ge1}\PR{M\left(\frac{y-w}{\mu}\right)=m}\Bigg]\\
   &&-a\int_z^y\!\!dw\,h_{n+1}(aw)\\
   &&\qquad\quad\left[\Theta_{m-n-1}\!\!\sum_{k=n+1}^{m-1}\PR{M\left(\frac{w-z}{\mu}\right)=m-k}\PR{M\left(\frac{y-w}{\mu}\right)=k}\right.\\
   &&\qquad\qquad\qquad\qquad-\PR{M\left(\frac{w-z}{\mu}\right)\ge1}\PR{M\left(\frac{y-w}{\mu}\right)=m}\Bigg]
\end{eqnarray*}
On the other hand, since for every $t>0$
\begin{equation*}
    \PR{M(t)\ge1}=1-\PR{M(t)=0}
\end{equation*}
we can also write
\begin{eqnarray*}
   p_{m,n} \!\!\! &=&\!\!a\int_0^y\!\!dw\,[h_n(aw)-h_{n+1}(aw)]\,\PR{M\left(\frac{y-z}{\mu}\wedge\frac{y-w}{\mu}\right)=m}\\
   &&-a\int_z^y\!\!dw\,[h_n(aw)-h_{n+1}(aw)]\,\PR{M\left(\frac{y-w}{\mu}\right)=m}\\
   &&+a\int_z^y\!\!dw\,h_n(aw)\sum_{k=n}^{m}\PR{M\left(\frac{w-z}{\mu}\right)=m-k}\PR{M\left(\frac{y-w}{\mu}\right)=k}\\
   &&-a\int_z^y\!\!dw\,h_{n+1}(aw)\!\sum_{k=n+1}^{m}\!\!\PR{M\left(\frac{w-z}{\mu}\right)=m-k}\!\PR{M\left(\frac{y-w}{\mu}\right)=k}\\
   &=&a\int_0^z\!\!dw\,[h_n(aw)-h_{n+1}(aw)]\,\PR{M\left(\frac{y-z}{\mu}\right)=m}\\
   &&\quad+a\int_z^y\!\!dw\Bigg\{h_{n+1}(aw)\PR{M\left(\frac{w-z}{\mu}\right)=m-n}\PR{M\left(\frac{y-w}{\mu}\right)=n}\\
   &&\qquad\qquad\quad+[h_n(aw)-h_{n+1}(aw)]\\
   &&\qquad\qquad\qquad\quad\sum_{k=n}^{m}\PR{M\left(\frac{w-z}{\mu}\right)=m-k}\PR{M\left(\frac{y-w}{\mu}\right)=k}\Bigg\}
\end{eqnarray*}
which when the Poisson probabilities are introduced goes into
\begin{eqnarray}
  p_{m,n}(y,z) &=& \pi_m(y-z)a\int_0^z\!\!dw\,[h_n(aw)-h_{n+1}(aw)]\nonumber\\
   &&\quad+a\int_z^y\!\!dw\, h_{n+1}(aw)\pi_{m-n}(w-z)\pi_n(y-w)\label{z>0m>n} \\
   &&\qquad\quad+a\int_z^y\!\!dw\,[h_n(aw)-h_{n+1}(aw)]\sum_{k=n}^m\pi_{m-k}(w-z)\pi_k(y-w)\nonumber\\
   &=&A_{m,n}(y,z)-A_{m,n+1}(y,z)+C_{m,n}(y,z)-C_{m,n+1}(y,z)\nonumber
\end{eqnarray}
namely the result of Proposition\myref{case1} once the
definitions\refeq{A1} and\refeq{C1} are taken into account

\subsection{Continuity in $z=0$}

We will now compare the values of $p_{m,n}(y,z)$ separately listed
in the Proposition\myref{case1} for $z<0$ and $z>0$, and we will
show that they connect with continuity in $z=0$, in the sense that
for every $y>0$
\begin{equation*}
    p_{m,n}(y,0^-)=p_{m,n}(y,0^+)
\end{equation*}
For $z<0$ (namely $a\mu s>\lambda t$) the results are given in the
Lemma\myref{lemma} and in formulas\refeq{z<0m=n} and \refeq{z<0m>n},
so that for $z\uparrow0^-$ we simply have
\begin{eqnarray}
  p_{m,n}(y,0^-) &=& 0\qquad\qquad\qquad\qquad\qquad\qquad\qquad\qquad\qquad\qquad n>m\ge0\qquad\label{pz<0n>m}\\
  p_{n,n}(y,0^-) &=& a\int_0^y\!\!dw\,h_n(aw)\,\pi_0(w)\,\pi_n(y-w)\qquad\qquad\qquad\;\;\; n=m\ge0\qquad\label{pz<0n=m}\\
  p_{m,n}(y,0^-) &=& a\int_0^y\!\!dw\Bigg\{h_{n+1}(aw)\pi_{m-n}(w)\pi_n(y-w)\ \qquad\qquad m>n\ge0\qquad\quad\label{pz<0n<m} \\
    &&\qquad+\big[h_n(aw)-h_{n+1}(aw)\big]\sum_{k=n}^{m}\pi_{m-k}(w)\pi_k(y-w)\Bigg\}\nonumber
\end{eqnarray}
On the other hand, when $z>0$ (namely $a\mu s<\lambda t$) the
results are given in formulas\refeq{z>0n>m},\refeq{z>0m=n} and
\refeq{z>0m>n}, where $z$ appears also as an integration limit so
that it is advisable to premise a few formal remarks. If the
integrand is a regular function $\xi(x)$ we of course have with
$y>z>0$
\begin{equation*}
    \lim_{z\downarrow0^+}\int_0^z\xi(x)\,dx=\int_0^{0^+}\xi(x)\,dx=0\qquad\qquad\lim_{z\downarrow0^+}\int_z^y\xi(x)\,dx=\int_0^y\xi(x)\,dx
\end{equation*}
When instead the integrand contains the distribution
$\delta(x)=\delta(x-0^+)$ concentrated on $x=0^+$, as indeed happens
in every first term of $h_n(x)$ which is
$\beta_0(n)\delta(x)=a^n\delta(x)$ (see also
Appendix\myref{notations}), we have
\begin{equation*}
    \lim_{z\downarrow0^+}\int_0^z\xi(x)\delta(x)\,dx=\int_0^{0^+}\xi(0)\delta(x)\,dx=\xi(0)\qquad\qquad\lim_{z\downarrow0^+}\int_z^y\xi(x)\delta(x)\,dx=0
\end{equation*}
As a consequence in the limit $z\downarrow0^+$, by neglecting the
vanishing terms, the following integrals on $[0,z]$ only retain the
contribution of the term $k=0$ in $h_n$
\begin{eqnarray*}
    \lim_{z\downarrow0^+}\int_0^zdx\,\xi(x)h_n(x)&=&\int_0^{0^+}dx\,\xi(x)\sum_{k=0}^n\beta_k(n)f_k(x)=\int_0^{0^+}\!\!dx\,\xi(x)\beta_0(n)f_0(x)\\
    &=&a^n\int_0^{0^+}\!\!dx\,\xi(0)\delta(x)=a^n\xi(0)
\end{eqnarray*}
while on the other hand for the integrals on $[z,y]$ we have
\begin{equation*}
    \lim_{z\downarrow0^+}\int_z^ydx\,\xi(x)h_n(x)=\int_{0^+}^ydx\,\xi(x)\sum_{k=0}^n\beta_k(n)f_k(x)=\int_0^ydx\,\xi(x)\sum_{k=1}^n\beta_k(n)f_k(x)
\end{equation*}
where now the term $k=0$ is apparently missing

Take first $p_{m,n}(y,z)$ for $z>0$ and $n>m\ge0$ as given
in\refeq{z>0n>m}: from the previous remarks, with $x=aw$ and
recalling also that $\pi_k(0)=\delta_{k,0}$ (so that
$\pi_{n-m}(0)=0$ because $n>m$), in the limit $z\downarrow0^+$ we
immediately have
\begin{eqnarray*}
    p_{m,n}(y,0^+)\!\!&=&\!\!\pi_m(y)\int_0^{0^+}\!\!dx\,\left[a^n-a^{n+1}+a^n\pi_{n-m}(0)-(a^n-a^{n+1})\sum_{k=0}^{n-m}\pi_k(0)\right]\delta(x) \\
    &=&\!\!\pi_m(y)\left[a^n-a^{n+1}-(a^n-a^{n+1})\sum_{k=0}^{n-m}\delta_{k,0}\right]=0
\end{eqnarray*}
which coincides with $p_{m,n}(y,0^-)$ in\refeq{pz<0n>m}. Then,
always with $z>0$, consider the case $n=m\ge0$ given
in\refeq{z>0m=n}: now (since $\pi_0(0)=1$) in the limit
$z\downarrow0^+$ we have that the first integral exactly compensates
the term $k=0$ missing in the sum of the second integral, so that
\begin{eqnarray*}
    p_{n,n}(y,0^+)&=&\pi_n(y)\int_0^{0^+}\!\!dx\,\left[a^n-a^{n+1}+a^{n+1}\pi_0(0)\right]\delta(x)\\
    &&\qquad\qquad\qquad+a\int_0^ydw\,\pi_0(w)\pi_n(y-w)\sum_{k=1}^n\beta_k(n)f_k(aw) \\
    &=&a^n\pi_n(y)+a\int_0^ydw\,\pi_0(w)\pi_n(y-w)\sum_{k=1}^n\beta_k(n)f_k(aw)
    \\
    &=&a\int_0^ydw\,\pi_0(w)\pi_n(y-w)\sum_{k=0}^n\beta_k(n)f_k(aw)\\
    &=&a\int_0^ydw\,\pi_0(w)\pi_n(y-w)h_n(aw)
\end{eqnarray*}
which again coincides with $p_{n,n}(y,0^-)$ in\refeq{pz<0n=m}.
Finally consider the case $z>0$ and $m>n\ge0$ which is given
in\refeq{z>0m>n}: here with the same remarks as before in the limit
$z\downarrow0^+$ we again have $p_{m,n}(y,0^-)$ of\refeq{pz<0n<m}
\begin{eqnarray*}
    p_{m,n}(y,0^+)\!\!&=&\!\!\pi_m(y)\int_0^{0^+}\!\!dx\,(a^n-a^{n+1})\,\delta(x) \\
    &&\qquad\qquad+a\int_{0^+}^ydw\Bigg\{h_{n+1}(aw)\pi_{m-n}(w)\pi_n(y-w)\\
    &&\qquad\qquad\qquad+[h_n(aw)-h_{n+1}(aw)]\sum_{k=n}^m\pi_{m-k}(w)\pi_k(y-w)\Bigg\}\\
    &=&\!\!a^n(1-a)\pi_m(y)+a\int_{0^+}^ydw\Bigg\{h_{n+1}(aw)\pi_{m-n}(w)\pi_n(y-w)\\
    &&\qquad\qquad\qquad+[h_n(aw)-h_{n+1}(aw)]\sum_{k=n}^m\pi_{m-k}(w)\pi_k(y-w)\Bigg\}\\
    &=&\!\!a\int_0^ydw\Bigg\{h_{n+1}(aw)\pi_{m-n}(w)\pi_n(y-w)\\
    &&\qquad\qquad\qquad+[h_n(aw)-h_{n+1}(aw)]\sum_{k=n}^m\pi_{m-k}(w)\pi_k(y-w)\Bigg\}
\end{eqnarray*}
because it is easy to see that the first term exactly compensates
the missing term in the second integral:
\begin{eqnarray*}
  \lefteqn{\int_0^{ay}dx\left[a^{n+1}\pi_{m-n}(0)\pi_n(y)+(a^n-a^{n+1})\sum_{k=n}^m\pi_{m-k}(0)\pi_k(y)\right]\delta(x)}\qquad\qquad\qquad\qquad \\
   &=&\int_0^{ay}dx\,a^n(1-a)\sum_{k=n}^m\delta_{m-k,0}\pi_k(y)=a^n(1-a)\pi_m(y)
\end{eqnarray*}

\section{A proof of Proposition\myref{case2}}\label{proof2}

To find an explicit elementary formula for $Q_{m,n}(y,z)$ (here
$z<0$, namely $a\mu s>\lambda t$) let us begin by supposing
$m>n\ge1$: in this case -- within our usual notations --
from\refeq{Q} we have
\begin{eqnarray*}
  Q_{m,n}(y,z) &=& \frac{e^{-(y-z)}}{m!}a\int_0^ydw\left[\beta_0(n)f_0(aw)+\sum_{\ell=1}^n\beta_\ell(n)f_\ell(aw)\right] \\
   &&\qquad\qquad\qquad\qquad\qquad\qquad\cdot\sum_{k=n}^m\binom{m}{k}(w-z)^{m-k}(y-w)^k\\
   &=&\frac{e^{-(y-z)}}{m!\,a^m}\int_0^{ay}dx\left[a^n\delta(x)+\sum_{\ell=1}^n\beta_\ell(n)f_\ell(x)\right] \\
   &&\qquad\qquad\qquad\qquad\qquad\qquad\cdot\sum_{k=n}^m\binom{m}{k}(x-az)^{m-k}(ay-x)^k\\
   &=&\frac{e^{-(y-z)}}{m!\,a^m}\sum_{k=n}^m\binom{m}{k}\Bigg[a^n\left(-az\right)^{m-k}\left(ay\right)^k\\
   &&\qquad\qquad+\sum_{\ell=1}^n\frac{\beta_\ell(n)}{(\ell-1)!}\int_0^{ay}dx\,e^{-x}x^{\ell-1}(x-az)^{m-k}(ay-x)^k\Bigg]
\end{eqnarray*}
On the other hand (see\mycite{grad} $3.383.1$) with $v=ay-x$ it is
\begin{eqnarray*}
  \lefteqn{\int_0^{ay}\!\!dx\,e^{-x}x^{\ell-1}(x-az)^{m-k}(ay-x)^k=e^{-ay}\!\!\int_0^{ay}\!\!dv\,e^v(ay-v)^{\ell-1}[a(y-z)-v]^{m-k}v^k}\qquad\qquad\quad\\
   &=&e^{-ay}\int_0^{ay}dv\,e^v(ay-v)^{\ell-1}v^k\sum_{i=0}^{m-k}\binom{m-k}{i}[a(y-z)]^i(-v)^{m-k-i}\\
   &=&e^{-ay}\sum_{i=0}^{m-k}\binom{m-k}{i}(-1)^{m-k-i}[a(y-z)]^i\int_0^{ay}dv\,e^v(ay-v)^{\ell-1}v^{m-i}\\
   &=&e^{-ay}\sum_{i=0}^{m-k}\binom{m-k}{i}(-1)^{m-k-i}[a(y-z)]^i(ay)^{\ell+m-i}\\
   &&\qquad\qquad\qquad\qquad\mathrm{B}(\ell,m-i+1)\,\Phi(m-i+1;m-i+\ell+1;ay)
\end{eqnarray*}
where $\Phi(\alpha;\beta;x)$ is a confluent hypergeometric function
(see\mycite{grad} $9.2$) and $\mathrm{B}(\alpha,\beta)$ is the beta
function (see\mycite{grad} $8.38$) with
\begin{equation*}
  \mathrm{B}(\alpha,\beta) =
  \frac{\Gamma(\alpha)\Gamma(\beta)}{\Gamma(\alpha+\beta)}=\frac{(\alpha-1)!(\beta-1)!}{(\alpha+\beta-1)!}
\end{equation*}
As a consequence we have
\begin{eqnarray*}
  Q_{m,n}(y,z) &=&\frac{e^{-(y-z)}}{m!\,a^m}\sum_{k=n}^m\binom{m}{k}\Bigg[a^n\left(-az\right)^{m-k}\left(ay\right)^k\\
   &&\qquad+\sum_{\ell=1}^n\beta_\ell(n)e^{-ay}\sum_{i=0}^{m-k}\binom{m-k}{i}(-1)^{m-k-i}[a(y-z)]^i(ay)^{\ell+m-i}\\
   &&\qquad\qquad\qquad\qquad\frac{(m-i)!}{(m-i+\ell)!}\,\Phi(m-i+1;m-i+\ell+1;ay)\Bigg]\\
   &=&\frac{e^{-(y-z)}}{m!\,a^m}\sum_{k=n}^m\binom{m}{k}\sum_{\ell=0}^n\beta_\ell(n)e^{-ay}\sum_{i=0}^{m-k}\binom{m-k}{i}[a(y-z)]^i(ay)^{\ell+m-i}\\
   &&\qquad\qquad(-1)^{m-k-i}\frac{(m-i)!}{(m-i+\ell)!}\,\Phi(m-i+1;m-i+\ell+1;ay)\Bigg]\\
   &=&\sum_{k=n}^m\sum_{i=0}^{m-k}\frac{(-1)^{m-k-i}}{a^{m-i}}\binom{m-i}{k}\sum_{\ell=0}^n\beta_\ell(n)\\
   &&\qquad\qquad\qquad\pi_i(y-z)\pi_{m-i+\ell}(ay)\Phi(m-i+1;m-i+\ell+1;ay)
\end{eqnarray*}
and with the change of the summation index $j=m-i$
\begin{eqnarray*}
  Q_{m,n}(y,z) &=&\sum_{k=n}^m\sum_{j=k}^m\frac{(-1)^{j-k}}{a^j}\binom{j}{k}\sum_{\ell=0}^n\beta_\ell(n)\\
   &&\qquad\qquad\qquad\pi_{m-j}(y-z)\pi_{j+\ell}(ay)\Phi(j+1;j+\ell+1;ay)
\end{eqnarray*}
as stated in the Proposition\myref{case2}. We have proved this
result by supposing $m>n\ge1$, but it is possible to check now by
direct calculation that also its extension to $m\ge n\ge0$ gives the
right results for $Q_{m,0}$ with $m\ge0$, and for $Q_{n,n}$ with
$n\ge1$. We have in fact from the definition\refeq{Q} that for $m\ge
n=0$ it is
\begin{eqnarray*}
  Q_{0,0}(y,z) &=& e^{-(y-z)}a\int_0^ydw\, \delta(aw)=e^{-(y-z)}=\pi_0(y-z) \\
  Q_{m,0}(y,z) &=& \frac{e^{-(y-z)}}{m!}\,a\int_0^ydw\,
   h_0(aw)\sum_{k=0}^{m}\binom{m}{k}(w-z)^{m-k}(y-w)^k
   \\
   &=&\frac{e^{-(y-z)}(y-z)^m}{m!}\,a\int_0^ydw\,\delta(aw)=\pi_m(y-z)
\end{eqnarray*}
while for $m=n\ge1$ we have with $v=a(y-w)$
\begin{eqnarray*}
  Q_{n,n}(y,z) &=& \frac{e^{-(y-z)}}{n!}\,a\int_0^ydw\left[a^n\delta(aw)+\sum_{k=1}^n\beta_k(n)f_k(aw)\right](y-w)^n \\
   &=&\frac{e^{-(y-z)}}{n!}\left[(ay)^n+\frac{e^{-ay}}{a^n}\sum_{k=1}^n\frac{\beta_k(n)}{(k-1)!}\int_0^{ay}e^vv^n(ay-v)^{k-1}dv\right]
\end{eqnarray*}
and since (see\mycite{grad} $3.383.1$)
\begin{equation*}
  \int_0^{ay}e^vv^n(ay-v)^{k-1}dv = (ay)^{k+n}\mathrm{B}(k,n+1)\Phi(n+1;n+k+1;ay)
\end{equation*}
we finally get
\begin{eqnarray*}
  Q_{n,n}(y,z) &=& \frac{e^{-(y-z)}}{n!}(ay)^n\left[1+\frac{e^{-ay}}{a^n}\sum_{k=1}^n\beta_k(n)\frac{n!(ay)^k}{(n+k)!}
            \Phi(n+1;n+k+1;ay)\right]
   \\
   &=&e^{-(y-z)}\pi_n(ay)\sum_{k=0}^n\frac{\beta_k(n)}{a^n}\frac{n!}{(n+k)!}(ay)^k
            \Phi(n+1;n+k+1;ay)
\end{eqnarray*}
It would be easy to see now that all these results can also be
derived as particular cases from our explicit expression\refeq{Q}
that can hence be considered as completely general

We proceed now to calculate for $z>0$ (namely for $\lambda t>a\mu
s$) the explicit elementary form of the terms $A, B$ and $C$ defined
in\refeq{A1},\refeq{B1} and\refeq{C1}. With $x=aw$, and by taking
into account\mycite{grad} $8.350.1$ and $8.352.1$, from\refeq{A1} we
first have for $m\ge0$ and $n\ge1$
\begin{eqnarray*}
  A_{m,n}(y,z) &=& \pi_m(y-z)\,a\int_0^zdw\left[\beta_0(n)f_0(aw)+\sum_{k=1}^n\beta_k(n)f_k(aw)\right] \\
   &=&\pi_m(y-z)\left[a^n+\sum_{k=1}^n\frac{\beta_k(n)}{(k-1)!}\int_0^{az}x^{k-1}e^{-x}dx\right]\\
%   &=&\pi_m(y-z)\left[a^n+\sum_{k=1}^n\frac{\beta_k(n)}{(k-1)!}\,\gamma(k,\lambda t-a\mu
%   s)\right]\\
   &=&\pi_m(y-z)\left\{a^n+\sum_{k=1}^n\beta_k(n)\left[1-\sum_{j=0}^{k-1}\pi_j(az)\right]\right\}\\
   &=&\pi_m(y-z)\sum_{k=0}^n\beta_k(n)\left[1+\pi_k(az)-\sum_{j=0}^k\pi_j(az)\right]
\end{eqnarray*}
Since however this result can be extended also to $n=0$ giving the
correct result
\begin{equation*}
    A_{m,0}(y,z)= \pi_m(y-z)\int_0^{az}\delta(x)\,dx=\pi_m(y-z)
\end{equation*}
it can be definitely used for every possible value of $n,m\ge0$ as
recalled in\refeq{A} of the Proposition\myref{case2}

Then, for $n\ge m\ge0$ and $x=aw$, by using\mycite{grad} $3.383.1$
we have from\refeq{B1}
\begin{eqnarray*}
B_{m,n}(y,z)&=& \pi_m(y-z)\,a\!\int_0^zdw\,h_{n+1}(aw)\sum_{k=0}^{n-m}\pi_k(z-w)\\
  &=& \pi_m(y-z)\sum_{k=0}^{n-m}\int_0^{az}dx\,
  h_{n+1}(x)\pi_k\left(z-\frac{x}{a}\right)\\
  &=&
  \pi_m(y-z)\sum_{k=0}^{n-m}\sum_{\ell=0}^{n+1}\beta_\ell(n+1)\int_0^{az}dx\,
  f_\ell(x)\pi_k\left(z-\frac{x}{a}\right)\\
  &=& \pi_m(y-z)\sum_{k=0}^{n-m}\Bigg[a^{n+1}\pi_k\left(z\right)\\
  &&\qquad\qquad\left.+\frac{e^{-z}}{a^kk!}\sum_{\ell=1}^{n+1}\frac{\beta_\ell(n+1)}{(\ell-1)!}\int_0^{az}dx\,x^{\ell-1}(az-x)^ke^{\frac{1-a}{a}x}\right]\\
  &=& \pi_m(y-z)\sum_{k=0}^{n-m}\pi_k\left(z\right)\Bigg[a^{n+1}\\
  &&\qquad\qquad\left.+\sum_{\ell=1}^{n+1}\beta_\ell(n+1)\frac{(az)^\ell
  k!}{(k+\ell)!}\,\Phi\left(\ell,k+\ell+1,(1-a)z\right)\right]\\
  &=&\pi_m(y-z)\sum_{k=0}^{n-m}\pi_k\left(z\right)\sum_{\ell=0}^{n+1}\beta_\ell(n+1)\frac{(az)^\ell
  k!}{(k+\ell)!}\,\Phi\left(\ell,k+\ell+1,(1-a)z\right)
\end{eqnarray*}
because $\Phi(0,k+1,x)=1$. This result also coincides with\refeq{B}
in the Proposition\myref{case2}

Finally, for $m\ge n\ge0$ and $0<z<y$, from\refeq{C1} we first
remark that $C_{m,0}(y,z)=0$ because $h_0(aw)=\delta(aw)$ is now
peaked outside the integration interval $[z,y]$. Then for $m\ge
n\ge1$ (the term of the sum with $\ell=0$ again vanishes for the
same reason as before) with $v=a(y-w)$, and $u=a(y-z)$ for short,
from\refeq{C1} we have
\begin{eqnarray*}
  C_{m,n}(y,z) &=& \int_0^udv\,h_n(ay-v)\sum_{k=n}^m\pi_{m-k}\left(\frac{u-v}{a}\right)\pi_k\left(\frac{v}{a}\right) \\
  &=&\sum_{\ell=1}^n\beta_\ell(n)\int_0^u\!\!dv\,f_\ell(ay-v)\sum_{k=n}^m\pi_{m-k}\left(\frac{u-v}{a}\right)\pi_k\left(\frac{v}{a}\right)\\
  &=& \frac{e^{-(y-z)}}{a^m}\sum_{\ell=1}^n\beta_\ell(n)\int_0^u\!\!dv\,f_\ell(ay-v)\sum_{k=n}^m\frac{(u-v)^{m-k}v^k}{(m-k)!k!}\\
  &=& \frac{e^{-(y-z)}}{a^mm!}\sum_{\ell=1}^n\beta_\ell(n)\sum_{k=n}^m\binom{m}{k}\int_0^u\!\!\!dv\,\frac{(ay-v)^{\ell-1}e^{-ay+v}}{(\ell-1)!}(u-v)^{m-k}v^k\\
  &=& \frac{e^{-(y-z)}e^{-ay}}{a^mm!}\sum_{\ell=1}^n\frac{\beta_\ell(n)}{(\ell-1)!}\sum_{k=n}^m\binom{m}{k}\int_0^u\!\!dv\,e^v(ay-v)^{\ell-1}(u-v)^{m-k}v^k
\end{eqnarray*}
and since it is (see\mycite{grad} $3.383.1$)
\begin{eqnarray*}
  \lefteqn{\int_0^udv\,e^v(ay-v)^{\ell-1}(u-v)^{m-k}v^k} \\
   &=& \sum_{j=0}^{\ell-1}\binom{\ell-1}{j}(ay)^j(-1)^{\ell-1-j}\int_0^u\!\!dv\,e^vv^{k+\ell-1-j}(u-v)^{m-k}\\
   &=& \sum_{j=0}^{\ell-1}\binom{\ell-1}{j}(ay)^j(-1)^{\ell-1-j}\\
   &&\qquad\quad B(m-k+1,k+\ell-j)u^{m+\ell-j}\Phi(k+\ell-j,m+\ell-j+1,u)\\
   &=& \sum_{j=0}^{\ell-1}\frac{(\ell-1)!(m-k)!(k+\ell-1-j)!}{j!(\ell-1-j)!(m+\ell-j)!}(-1)^{\ell-1-j}\\
   &&\qquad\qquad\qquad\qquad\qquad\qquad\quad (ay)^ju^{m+\ell-j}\Phi(k+\ell-j,m+\ell-j+1,u)
\end{eqnarray*}
we can write
\begin{eqnarray*}
  C_{m,n}(y,z) &=& \frac{e^{-(y-z)}e^{-ay}}{a^m}\sum_{\ell=1}^n\beta_\ell(n)\sum_{k=n}^m\sum_{j=0}^{\ell-1}\binom{k+\ell-j-1}{k}\\
  &&\qquad\quad\frac{(-1)^{\ell-1-j}(ay)^ju^{m+\ell-j}}{j!(m+\ell-j)!}\,\Phi(k+\ell-j,m+\ell-j+1,u) \\
  &=&\frac{e^{-(1-a)(y-z)}}{a^m}\sum_{\ell=1}^n\beta_\ell(n)\sum_{k=n}^m\sum_{j=0}^{\ell-1}(-1)^{\ell-1-j}\binom{k+\ell-j-1}{k}\\
  &&\qquad \pi_j(ay)\pi_{m+\ell-j}(a(y-z))\Phi(k+\ell-j,m+\ell-j+1,a(y-z))
\end{eqnarray*}
as stated in\refeq{C} of the Proposition\myref{case2}

As for the confluent hypergeometric functions appearing in all these
results they in fact boil down to finite combinations of elementary
functions as stated in the final formulas of the
Proposition\myref{case2}. The case $n=0$ immediately follows indeed
from the definitions (see\mycite{grad} $9.210.1$). On the other
hand, when $1\le\alpha\le\beta$, with the change of variables
$y=x-z$ and by taking sequentially into account\mycite{grad}
$3.383.1$, $8.350.1$ and $8.352.1$ we have
\begin{eqnarray*}
  \Phi(\alpha,\beta+1,x) &=& \frac{1}{B(\beta-\alpha+1,\alpha)x^\beta}\int_0^xe^zz^{\alpha-1}(x-z)^{\beta-\alpha}dz \\
   &=&\frac{\beta!\,e^x}{(\beta-\alpha)!(\alpha-1)!\,x^\beta}\int_0^xe^{-y}(x-y)^{\alpha-1}y^{\beta-\alpha}dy
   \\
   &=&\frac{\beta!\,e^x}{(m-\alpha)!(\alpha-1)!\,x^\beta}\sum_{\gamma=0}^{\alpha-1}\binom{\alpha-1}{\gamma}x^\gamma(-1)^{\alpha-\gamma-1}\!\!\int_0^x\!\!y^{\beta-\gamma-1}e^{-y}dy\\
%   &=&e^x\sum_{k=0}^{\alpha-1}\frac{(-1)^{\alpha-k-1}}{(\beta-\alpha)!(\alpha-k-1)!}\,\frac{\pi_k(x)}{\pi_\beta(x)}\,\gamma(\beta-k,x)\\
   &=&e^x\sum_{\gamma=0}^{\alpha-1}\frac{(-1)^{\alpha-\gamma-1}(\beta-\gamma-1)!}{(\beta-\alpha)!(\alpha-\gamma-1)!}\,\frac{\pi_\gamma(x)}{\pi_\beta(x)}\left(1-\sum_{\eta=0}^{\alpha-\gamma-1}\pi_\eta(x)\right)\\
   &=&e^x\sum_{\gamma=0}^{\alpha-1}(-1)^{\alpha-\gamma-1}\binom{\beta-\gamma-1}{\beta-\alpha}\frac{\pi_\gamma(x)}{\pi_\beta(x)}\left(1-\!\!\sum_{\eta=0}^{\alpha-\gamma-1}\pi_\eta(x)\!\right)
\end{eqnarray*}
as stated in the proposition

\end{appendix}

\end{document}